\documentclass[reqno,12pt]{amsart}
\usepackage{amscd,amsmath,amssymb}
\theoremstyle{plain}
\newtheorem{thm}{Theorem}[section]
\theoremstyle{plain}
\newtheorem{cor}[thm]{Corollary} 
\newtheorem{lemma}[thm]{Lemma} 
\newtheorem{prop}[thm]{Proposition}
\newtheorem{defi}[thm]{Definition}

\newtheorem{example}[thm]{Example}

\newtheorem{remark}[thm]{Remark}

\newcommand\Ad{{\operatorname{Ad}}}
\newcommand\Alg{{\operatorname{Alg}}}
\newcommand\Span{{\operatorname{Span}}}
\newcommand\card{{\operatorname{card}}}

\begin{document}
\title{ Reduced Free Products of Finite Dimensional $C^*$-Algebras } \par
\author{ Nikolay A. Ivanov }
\date{\today}

\address{\hskip-\parindent
Nikolay Ivanov  \\
Department of Mathematics \\
Texas A\&M University \\
College Station TX 77843-3368, USA}
\email{nivanov@math.tamu.edu}

\begin{abstract}
We find a necessary and sufficient conditions for the simplicity and uniqueness of trace for reduced free 
products of finite families of finite dimensional $C^*$-algebras with specified traces on them.
\end{abstract}

\maketitle

\section{Introduction and Definitions}

The notion of reduced free product of a family of $C^*$-algebras with specified states on them was introduced independently by Avitzour
 (\cite{A82}) and Voiculescu (\cite{V85}). We will recall this notion and some of its properties here.
\par

\begin{defi}
  The couple $(A,\phi)$, where $A$ is a unital $C^*$-algebra and $\phi$ a state is called a $C^*$-noncommutative probability space or $C^*$-NCPS.
\end{defi}

\par

\begin{defi}
  Let $(A,\phi)$ be a $C^*$-NCPS  and $\{ A_i | i \in I \}$ be a family of $C^*$-subalgebras of $A$, s.t. $1_A \in A_i$, $\forall i\in I$, where $I$ is an index set. We say that the family $\{ A_i |i \in I \}$ is free if $\phi(a_1...a_n)=0$, whenever $a_j \in A_{i_j} $ with $i_1\neq i_2\neq ... \neq i_n$ and $\phi(a_j)=0$, $\forall j \in \{ 1,...n \}$.
  A family of subsets $\{ S_i | i \in I \}$ $\subset$ $A$ is $*$-free if 
  $\{ C^*(S_i \cup \{ 1_A \} ) | i \in I \}$ is free.
\end{defi}

Let $\{ (A_i,\phi_i) | i \in I \}$ be a family of $C^*$-NCPS such that the GNS representations of $A_i$ associated to $\phi_i$ are all faithful. Then there is a unique $C^*$-NCPS $(A,\phi) \overset{def}{=} \underset{i \in I}{*} (A_i,\phi_i)$ with unital embeddings $A_i \hookrightarrow A$, s.t.
\\
(1) $\phi|_{A_i}=\phi_i$
\\
(2) the family $\{ A_i | i \in I \}$ is free in $(A,\phi)$
\\
(3) $A$ is the $C^*$-algebra generated by $\underset{i \in I}{\bigcup}A_i$
\\
(4) the GNS representation of $A$ associated to $\phi$ is faithful.
\\
And also:
\\
(5) If $\phi_i$ are all traces then $\phi$ is a trace too (\cite{V85}).
\\
(6) If $\phi_i$ are all faithful then $\phi$ is faithful too (\cite{D98}).

\par
In the above situation $A$ is called the reduced free product algebra and $\phi$ is called the free product state. Also the construction of the
reduced free product is based on defining a free product Hilbert space, which turns out to be $\mathfrak{H}_A$ - the GNS Hilbert space for 
$A$, associated to $\phi$.
\par

\begin{example}
If $\{ G_i | i \in I \}$ is a family of discrete groups and $C^*_r(G_i)$ are the reduced group $C^*$-algebras, corresponding to the left 
regular representations of $G_i$ on $l^2(G_i)$ respectively, and if $\tau_i$ are the canonical traces on $C^*_r(G_i)$, $i \in I$, then 
we have $\underset{i \in I}{*} (C^*_r(G_i), \tau_i)=(C^*_r(\underset{i \in I}{*}G_i), \tau)$, where $\tau$ is the canonical trace on the group $C^*$-algebra $C^*_r(\underset{i \in I}{*} G_i)$.
\end{example}

Reduced free products satisfy the following property:

\begin{lemma}[\cite{DR98}]

Let $I$ be an index set and let $(A_i,\phi_i)$ be a $C^*$-NCPS ($i \in I$), where each $\phi_i$ is faithful. Let $(B,\psi)$ be a $C^*$-NCPS
with $\psi$ faithful. Let 

\begin{center}

$(A,\phi) = \underset{i\in I}{*} (A_i,\phi_i)$.

\end{center}

Given unital $*$-homomorphisms, $\pi_i : A_i \rightarrow B$, such that $\psi \circ \pi_i = \phi_i$ and $\{ \pi_i(A_i) \}_{i\in I}$ is free
in $(B, \psi)$, there is a $*$-homomorphism, $\pi : A \rightarrow B$ such that $\pi|_{A_i} = \pi$ and $\psi \circ \pi = \phi$.

\end{lemma}

\par

From now on we will be concerned only with $C^*$-algebras equipped with tracial states.

\par

The study of simplicity and uniqueness of trace for reduced free
products of $C^*$-algebras, one can say, started with the paper of
Powers \cite{P75}. In this paper Powers proved that the reduced 
$C^*$-algebra of the free group on two generators $F_2$ is simple and
has a unique trace - the canonical one. In \cite{C79} Choi showed the same for the
"Choi algebra" $C_r^*(\mathbb{Z}_2 * \mathbb{Z}_3)$ and then
Paschke and Salinas in \cite{PS79} generalized the result to the case of $C_r^*(G_1 *
G_2)$, where $G_1, G_2$ are discrete groups, such that $G_1$ has
at least two and $G_2$ at least three elements. After that Avitzour in \cite{A82} gave a sufficient condition for simplicity and uniqueness of
trace for reduced free products of $C^*$-algebras, generalizing the
previous results. He proved:

\begin{thm}[\cite{A82}]

Let 

\begin{equation*}
(\mathfrak{A}, \tau) = (A, \tau_A) * (B, \tau_B),
\end{equation*}
where $\tau_A$ and $\tau_B$ are traces and $(A,\tau_A)$ and $(B,\tau_B)$ have faithful GNS representations. Suppose
that there are unitaries $u,v \in A$ and $w \in B$, such that 
$\tau_A(u) = \tau_A(v) = \tau_A(u^* v) = 0$ and $\tau_B(w) = 0$. Then
$\mathfrak{A}$ is simple and has a unique trace $\tau$.

\end{thm} 

{\em Note:} It is clear that $uw$ satisfies $\tau((uw)^n) = 0$, $\forall n \in \mathbb{Z} \backslash \{ 0 \}$. Unitaries with this property we
define below.

\section{Statement of the Main Result and Preliminaries}

We adopt the following notation:
\\
If $A_0$, ... , $A_n$ are unital $C^*$-algebras equipped with traces $\tau_0$, ... , $\tau_n$ respectively, then $A=\underset{\alpha_0}{\overset{p_0}{A_0}} \bigoplus
\underset{\alpha_1}{\overset{p_1}{A_1}} \bigoplus ... \bigoplus \underset{\alpha_n}{\overset{p_n}{A_n}}$ will mean that the $C^*$-algebra
$A$ is isomorphic to the direct sum of $A_0$, ... , $A_n$, and is such that $A_i$ are supported on the projections $p_i$. Also $A$ comes with a
trace (let's call it $\tau$) given by the formula $\tau=\alpha_0\tau_0 + \alpha_1\tau_1 + ... + \alpha_n\tau_n$. Here of course $\alpha_0$, $\alpha_1$, ... ,
$\alpha_n > 0$ and $\alpha_0 + \alpha_1 + ... + \alpha_n = 1$.

\begin{defi}

If $(A,\tau)$ is a $C^*$-NCPS and $u\in A$ is a unitary with $\tau(u^n)=0$, $\forall n \in \mathbb{Z} \backslash \{ 0 \}$, then we call $u$ a 
Haar unitary.

\par

If $1_A \in B \subset A$ is a unital abelian $C^*$-subalgebra of $A$ we call $B$ a diffuse abelian $C^*$-subalgebra of $A$ if $\tau|_B$ is
given by an atomless measure on the spectrum of $B$. We also call $B$ a unital diffuse abelian $C^*$-algebra.

\end{defi}

From Proposition 4.1(i), Proposition 4.3 of \cite{DHR97} we can conclude the following:

\begin{prop}

If $(B,\tau)$ is a $C^*$-NCPS with $B$-abelian, then $B$ is diffuse abelian if and only if $B$ contains a Haar unitary.

\end{prop}

$C^*$-algebras of the form $(\underset{\alpha}{\overset{p}{\mathbb{C}}} \bigoplus
\underset{1-\alpha}{\overset{1-p}{\mathbb{C}}})*(\underset{\beta}{\overset{q}{\mathbb{C}}}\bigoplus
\underset{1-\beta}{\overset{1-q}{\mathbb{C}}})$ have been described explicitly in \cite{ABH91} (see also \cite{D99LN}):

\begin{thm}

Let $1 > \alpha \geqq \beta \geqq \frac{1}{2}$ and let

\begin{center}

$( A,\tau ) = ( \underset{\alpha }{\overset{p}{\mathbb{C}}} \oplus
\underset{1-\alpha }{\overset{1-p}{\mathbb{C}}} ) * ( \underset{\beta }{\overset{q}{\mathbb{C}}}\oplus
\underset{1-\beta}{\overset{1-q}{\mathbb{C}}} ) $.

\end{center}

If $\alpha > \beta$ then

\begin{equation*} 
A=\underset{\alpha -\beta }{\overset{p\wedge (1-q)}{\mathbb{C}}}\oplus
C([a,b], M_2(\mathbb{C}))\oplus \underset{\alpha + \beta -1}{\overset{p\wedge q}{\mathbb{C}}} ,
\end{equation*}
for some $0 < a < b < 1$. Furthermore, in the above picture

\begin{center}
$p=1 \oplus  \begin{pmatrix} 1 & 0 \\ 0 & 0 \end{pmatrix} \oplus 1 ,$
\end{center}

\begin{equation*}
q=0\oplus  \begin{pmatrix} t & \sqrt{t(1-t)} \\ \sqrt{t(1-t)} & 1-t \end{pmatrix} \oplus 1 ,
\end{equation*}
and the faithful trace $\tau$ is given by the indicated weights on the projections $p\wedge (1-q)$ and $p\wedge q$, together with an
atomless measure, whose support is $[a,b]$.

\par

If $\alpha =\beta > \frac{1}{2}$ then

\begin{equation*}
A=\{\ f:[0,b]\rightarrow M_2(\mathbb{C}) |\ f\ is\ continuous\ and\ f(0)\ is\ diagonal\ \} \oplus 
\underset{\alpha + \beta -1}{\overset{p\wedge q}{\mathbb{C}}},
\end{equation*}
for some $0 < b < 1$. Furthermore, in the above picture

\begin{center}
$p= \begin{pmatrix} 1 & 0 \\ 0 & 0 \end{pmatrix} \oplus 1,$
\end{center}

\begin{equation*}
q= \begin{pmatrix} {t} & {\sqrt{t(1-t)}} \\ {\sqrt{t(1-t)}} & {1-t} \end{pmatrix} \oplus 1,
\end{equation*}
and the faithful trace $\tau$ is given by the indicated weight on the projection $p\wedge q$, together with an atomless measure on $[0,b]$.

\par

If $\alpha = \beta = \frac{1}{2}$ then

\begin{equation*}
A=\{\ f:[0,1]\rightarrow M_2(\mathbb{C}) |\ f\ is\ continuous\ and\ f(0)\ and\ f(1)\ are\ diagonal\ \}.
\end{equation*}

Furthermore in the above picture

\begin{center}
$p= \begin{pmatrix} 1 & 0 \\ 0 & 0 \end{pmatrix} ,$
\end{center}

\begin{equation*}
q= \begin {pmatrix} t & \sqrt{t(1-t)} \\ \sqrt{t(1-t)} & 1-t \end{pmatrix} ,
\end{equation*}
and the faithful trace $\tau$ is given by an atomless measure, whose support is $[0,1]$.

\end{thm}

The question of describing the reduced free product of a finite family of finite dimensional abelian $C^*$-algebras was studied by Dykema in
\cite{D99}. He proved the following theorem:

\begin{thm}[\cite{D99}]

Let 

\begin{equation*}
(\mathfrak{A},\phi )=(\underset{\alpha_0}{\overset{p_0}{A_0}} \oplus \underset{\alpha_1}{\overset{p_1}{\mathbb{C}}} \oplus ...
\oplus \underset{\alpha_n}{\overset{p_n}{\mathbb{C}}})*(\underset{\beta_0}{\overset{q_0}{B_0}} \oplus 
\underset{\beta_1}{\overset{q_1}{\mathbb{C}}} \oplus ... 
\oplus \underset{\beta_m}{\overset{q_m}{\mathbb{C}}}),
\end{equation*}
where $\alpha_0 \geq 0$ and $\beta_0 \geq 0$ and $A_0$ and $B_0$ are equipped with traces $\phi(p_0)^{-1} \phi|_{A_0}$, $\phi(q_0)^{-1} 
\phi|_{B_0}$ and $A_0$ and $B_0$ have diffuse abelian $C^*$-subalgebras, and where $n \geq 1$, $m \geq 1$ (if $\alpha_0 = 0$ or $\beta_0 =
0$, or both, then, of course, we don't impose any conditions on $A_0$ or $B_0$, or both respectively). 
Suppose also that $\dim(A) \geq 2$, $\dim(B) \geq 2$, and $\dim(A) + \dim(B) \geq 5$.

\par

Then

\begin{equation*}
\mathfrak{A} = \overset{r_0}{\mathfrak{A}_0} \oplus \underset{(i',j)\in L_+}{\bigoplus} 
\underset{\alpha_i + \beta_i -1}{\overset{p_i \wedge q_j}{\mathbb{C}}},
\end{equation*}
where $L_+ = \{ (i,j)| 1 \leq i \leq n$, $1 \leq j \leq m$ and $\alpha_i + \beta_j > 1 \}$, and where 
$\mathfrak{A}_0$ has a unital, diffuse abelian sublagebra supported on $r_0 p_1$ and another one supported on $r_0 q_1$.

\par

Let $ L_0 = \{(i,j)| 1 \leq i \leq n$, $1 \leq j \leq m$ and $\alpha_i + \beta_j = 1 \} .$ 

\par

If $L_0$ is empty then $\mathfrak{A}_0$ is simple and $\phi(r_0)^{-1} \phi|_{\mathfrak{A}_{0}}$ is the unique trace on 
$\mathfrak{A}_0.$

\par

If $L_0$ is not empty, then for each $(i,j) \in L_0$ there is a $*$-homomorphism $\pi_{(i,j)}: \mathfrak{A}_0 \rightarrow \mathbb{C}$ such
that $\pi_{(i,j)}(r_0 p_i) = 1 = \pi_{(i,j)}(r_0 q_j).$ Then: \\ 

(1)  $\mathfrak{A}_{00} \overset{def}{=} \underset{(i,j)\in L_0}{\bigcap} \ker (\pi_{(i,j)})$ \\ 

is simple and nonunital, and $\phi(r_0)^{-1} \phi|_{\mathfrak{A}_{00}}$ is the unique trace on 
$\mathfrak{A}_{00}.$ \\ 

(2)  For each $i\in \{1,...n \}, \ r_0 p_i$ is full in $\mathfrak{A}_0 \cap 
\underset{i' \neq i}{\underset{(i',j) \in L_0}{\bigcap}} \ker
(\pi_{(i',j)}).$ \\ 

(3)  For each $j \in \{ 1, ... , m \}, \ r_0 q_j$ is full in $\mathfrak{A}_{0} \cap 
\underset{j' \neq j}{\underset{(i,j') \in
L_0}{\bigcap}} \ker (\pi_{(i,j')}).$

\end{thm}

One can define von Neumann algebra free products, similarly to reduced free products of $C^*$-algebras.
We will denote by $\mathbb{M}_n$ the $C^*$-algebra (von Neumann algebra) of $n \times n$ matrices 
with complex coefficients.
\par
Dykema studied the case of von Neumann algebra free products of finite dimensional (von Neumann) 
algebras:

\begin{thm}[\cite{D93}]

Let 

\begin{equation*}
A = \underset{\alpha_0}{\overset{p_0}{L(F_s)}} \oplus \underset{\alpha_1}{\overset{p_1}{\mathbb{M}_{n_1}}} \oplus ... \oplus
\underset{\alpha_k}{\overset{p_k}{\mathbb{M}_{n_k}}}
\end{equation*}
and

\begin{equation*}
B = \underset{\beta_0}{\overset{q_0}{L(F_r)}} \oplus \underset{\beta_1}{\overset{q_1}{\mathbb{M}_{m_1}}} \oplus ... \oplus
\underset{\beta_l}{\overset{q_l}{\mathbb{M}_{m_l}}},
\end{equation*}
where $L(F_s), L(F_r)$ are interpolated free group factors, $\alpha_0, \beta_0 \geq 0$, 
and where $\dim(A) \geq 2$, $\dim(B) \geq 2$ and $\dim(A) + \dim(B)\geq 5$. 
Then for the von Neumann algebra free product we have:

\begin{equation*}
A*B = L(F_t) \oplus \underset{(i,j) \in L_+}{\bigoplus}
\underset{\gamma_{ij}}{\overset{f_{ij}}{\mathbb{M}_{N(i,j)}}},
\end{equation*}
where $L_+ = \{(i,j) | 1 \leq i \leq k, 1 \leq j \leq l, (\frac{\alpha_i}{n_i^2}) + (\frac{\beta_j}{m_j^2}) 
> 1 \}$, $N(i,j) = max(n_i, m_j)$, $\gamma_{ij} = N(i,j)^2 \cdot (\frac{\alpha_i}{n_i^2} + 
\frac{\beta_j}{m_j^2} - 1)$, and $f_{ij} \leq p_i \wedge q_j$.

\end{thm}

{\em Note:} $t$ can be determined from the other data, which makes sense only if the interpolated free group factors are all different. We
will use only the fact that $L(F_t)$ is a factor. For definitions and properties of interpolated free group factors see \cite{Ra94} and \cite{D94}.
\par
In this paper we will extend the result of Theorem 2.4 to the case of reduced free products of finite dimensional $C^*$-algebras with specified traces on them. We will prove:

\begin{thm}

Let

\begin{equation*}
(\mathfrak{A},\phi )=(\underset{\alpha_0}{\overset{p_0}{A_0}} \oplus \underset{\alpha_1}{\overset{p_1}{\mathbb{M}_{n_1}}} \oplus ...
\oplus \underset{\alpha_k}{\overset{p_k}{\mathbb{M}_{n_k}}})*(\underset{\beta_0}{\overset{q_0}{B_0}} \oplus 
\underset{\beta_1}{\overset{q_1}{\mathbb{M}_{m_1}}} \oplus ... \oplus \underset{\beta_l}{\overset{q_l}
{\mathbb{M}_{m_l}}}),
\end{equation*}
where $\alpha_0, \beta_0 \geq 0$, $\alpha_i > 0$, for $i=1,..,k$ and $\beta_j > 0$, for $j=1,...,l$, and where $\phi(p_0)^{-1} \phi|_{A_0}$ and
$\phi(q_0)^{-1} \phi|_{B_0}$ are traces on $A_0$ and $B_0$ respectivelly. Suppose that $\dim(A) \geq 2$, $\dim(B)
\geq 2$, $\dim(A) + \dim(B) \geq 5$, and that both $A_0$ and $B_0$ contain unital, diffuse abelian 
$C^*$-subalgebras (if $\alpha_0 > 0$, respectivelly $\beta_0 > 0$). Then

\begin{equation*}
\mathfrak{A}= \underset{\gamma}{\overset{f}{\mathfrak{A}_0}} \oplus \underset{(i,j)\in L_+}{\bigoplus}
\underset{\gamma_{ij}}{\overset{f_{ij}}{\mathbb{M}_{N(i,j)}}},	
\end{equation*}
where $L_+ = \{ (i,j)| \frac{\alpha_i}{n_i^2} + \frac{\beta_j}{m_j^2} > 1 \}$, $N(i,j) = max(n_i,m_j)$, $\gamma_{ij} =
N(i,j)^2(\frac{\alpha_i}{n_i^2} + \frac{\beta_j}{m_j^2} -1)$, $f_{ij} \leq p_i \wedge q_j$. There is a unital, diffuse abelian 
$C^*$-subalgebra of $\mathfrak{A}_0$, supported on $f p_1$ and another one, supported on $f q_1$.

\par

If $L_0 = \{ (i,j)| \frac{\alpha_i}{n_i^2} + \frac{\beta_j}{m_j^2} = 1 \},$ is empty, then $\mathfrak{A}_0$ is simple with a unique trace. If $L_0$ is not empty, then $\forall (i,j) \in L_0 ,\ \exists \pi_{(i,j)} :
\mathfrak{A}_{0} \rightarrow \mathbb{M}_{N(i,j)}$ a unital $*$-homomorphism, such that $\pi_{(i,j)}(f p_i) = \pi_{(i,j)}(f q_j) = 1$. Then: \\ 
(1) $\mathfrak{A}_{00} \overset{def}{=} \underset{(i,j) 
\in L_0}{\bigcap} \ker (\pi_{(i,j)})$ is simple and nonunital, and has a unique trace $\phi(f )^{-1} \phi |_{\mathfrak{A}_{00}}$. \\ 
(2) For each $i \in \{ 1, ..., k \}$, $f p_i$ is full in $\mathfrak{A}_0 \cap 
\underset{i' \neq i}{\underset{(i',j) \in L_0}{\bigcap}} \ker(\pi_{(i',j)})$. \\ 
(3) For each $j \in \{ 1, ..., l \}$, $f q_j$ is full in $\mathfrak{A}_0 \cap 
\underset{j' \neq j}{\underset{(i,j') \in L_0}{\bigcap}} \ker(\pi_{(i,j')})$. 

\end{thm}

\section{Beginning of the Proof - A Special Case}

In order to prove this theorem we will start with a simpler case. We will study first the $C^*$-algebras of the form 
$(A,\tau) \overset{def}{=} $ $ ( \underset{\alpha_1}{\overset{p_1}{\mathbb{C}}} \oplus ... \oplus \underset{\alpha_m}{\overset{p_m}
{\mathbb{C}}})*(\mathbb{M}_n, tr_n)$ with $0 < \alpha_1 \leq ... \leq \alpha_m$. We chose a set of matrix units for $\mathbb{M}_n$ and denote them by $\{ e_{ij}|i,j \in \{1,...n \} 
\} $ as usual. Let's take the (trace zero) permutation unitary
$$ u \overset{def}{=} \begin{pmatrix} 0 & 1 & ... & 0 \\ . & . & . & . \\ 0 & 0 & ... & 1 \\  1 & 0 & ... & 0 \end{pmatrix} \in \mathbb{M}_n.$$ \\ 

We see that $\Ad(u)(e_{11}) = u e_{11} u^* = e_{nn}$ and for $2 \leq i \leq n$, $\Ad(u)(e_{ii}) = u e_{ii} u^* = e_{(i-1) (i-1)}$. 
\par
It's clear that $$A = C^*(\{p_1, ..., p_m \}, \{ e_{ii}\}_{i=1}^n, u).$$ Then it is also clear that  
 $$A = C^*(\{ u^ip_1u^{-i} , ... , u^ip_mu^{-i} \}_{i=0}^{n-1}, \{ e_{ij} \}_{i=1}^{n}, u).$$ We want to show that the family 
 $$\{ \{ \mathbb{C} \cdot u^ip_1u^{-i}  \oplus ,..., \oplus \mathbb{C} \cdot u^ip_mu^{-i}  \}_{i=0}^{n-1},\ \{ \mathbb{C} \cdot e_{11}  \oplus ... \oplus \mathbb{C} \cdot e_{nn}  \} \}$$ is free.

We will prove something more general. We denote 
$$B \overset{def}{=} C^*( \{ u^kp_1u^{-k}, ... , u^kp_mu^{-k} \}_{k=0}^{n-1}, \{ e_{11}, ... ,e_{nn} \} ).$$ 
Let $l$ be an integer and $l|n$,  $1 < l < n$ (if such $l$ exists). Let $$E \overset{def}{=} C^*( \{ \{ u^kp_1u^{-k}, ... , u^kp_mu^{-k} \}_{k=0}^{l-1}, \{ e_{11} , ... , e_{nn} \}, \{ u^l, u^{2l}, ... , u^{n-l} \} \} ).$$ It's easy to see that
$$C^* ( \{ e_{11}, ... , e_{nn} \}, \{ u^l, u^{2l}, ... , u^{n-l} \} )= 
\underbrace{\mathbb{M}_{ \frac{n}{l} } \oplus ... \oplus \mathbb{M}_{ \frac{n}{l} }}_{l-times} \subset \mathbb{M}_n.$$ 
We will adopt the following notation from \cite {D99LN}:
\par
Let $(D, \varphi)$ be a $C^*$-NCPS and $1_D \in D_1, ..., D_k \subset D$ be a family of unital 
$C^*$-subalgebras of $D$, having a common unit $1_D$. We denote by $D^{\circ} \overset{def}{=} \{ d\in D |
\varphi(d)=0 \}$ (analoguously for $D_1$, ..., etc). We denote by $\Lambda^{\circ}(D_1^{\circ}, D_2^{\circ}
, ..., D_k^{\circ})$ 
the set of all words of the form $d_1 d_2 \cdots d_j$ and of nonzero length, where $d_t \in D_{i_t}^{\circ}$, for 
some $1 \leq i_t \leq k$ and $i_t \neq i_{t+1}$ for any $1 \leq t \leq j-1$. \\ 
\par
We have the following 

\begin{lemma}

If everything is as above, then:

(i)  The family $\{ \{ u^kp_1u^{-k} , ... , u^kp_mu^{-k} \}_{k=0}^{n-1},$ $\{ e_{11} , ... ,e_{nn} \} \}$ 
is free in $(A,\tau)$. And more generally if 
$$\omega \in \Lambda^{\circ}( C^*(p_1, ..., p_m)^{\circ}, ..., 
C^*(u^{n-1}p_1u^{1-n}, ..., u^{n-1}p_mu^{1-n})^{\circ}, C^*(e_{11}, ..., e_{nn})^{\circ}),$$ 
then $\tau(\omega u^r)=0$ for all $0 \leq r \leq n-1$.

(ii)  The family $\{ \{ u^kp_1u^{-k} , ... , uu^kp_mu^{-k} \}_{k=0}^{l-1},$ $\{ e_{11} , ... , e_{nn} , u^l, 
u^{2l}, ... u^{n-l} \} \}$ is free in $(A,\tau)$. And more generally if 
$$\omega \in \Lambda^{\circ}(C^*(p_1,..., p_m)^{\circ},..., 
C^*(u^{l-1}p_1u^{1-l},..., u^{l-1}p_mu^{1-l})^{\circ}, 
C^*( e_{11}, ..., e_{nn}, u^l,..., u^{n-l})^{\circ}),$$ 
then $\tau(\omega u^r)=0$ for all $0 \leq r \leq l-1$.

\end{lemma}

\begin{proof}

Each letter $\alpha \in C^*( \{ u^kp_1u^{-k}, ... , u^kp_mu^{-k} \})$ with $\tau(\alpha) = 0$ can be represented as $\alpha = u^k \alpha' u^{-k}$ with $\tau(\alpha') = 0$, and $\alpha' \in C^*( \{ p_1, ..., p_m \} )$. 
\par

Case (i): \\
\par
Each $$\omega \in \Lambda^{\circ}( C^*(p_1, ..., p_m)^{\circ}, ..., 
C^*(u^{n-1}p_1u^{1-n}, ..., u^{n-1}p_mu^{1-n})^{\circ}, C^*(e_{11}, ..., e_{nn})^{\circ})$$ is of one of the
four following types:

\begin{equation}
\omega =  \alpha_{11} \alpha_{12} \cdots \alpha_{1i_1} \beta_1 \alpha_{21} \cdots \alpha_{2i_2} \beta_2
\alpha_{31} \cdots \alpha_{t-1i_{t-1}} \beta_{t-1} \alpha_{t1} \cdots \alpha_{ti_t},
\end{equation}

\begin{equation}
\omega =  \beta_1 \alpha_{21} \cdots \alpha_{2i_2} \beta_2
\alpha_{31} \cdots \alpha_{t-1i_{t-1}} \beta_{t-1} \alpha_{t1} \cdots \alpha_{ti_t},
\end{equation}

\begin{equation}
\omega =  \beta_1 \alpha_{21} \cdots \alpha_{2i_2} \beta_2
\alpha_{31} \cdots \alpha_{t-1i_{t-1}} \beta_{t-1},
\end{equation}

\begin{equation}
\omega =  \alpha_{11} \alpha_{12} \cdots \alpha_{1i_1} \beta_1 \alpha_{21} \cdots \alpha_{2i_2} \beta_2
\alpha_{31} \cdots \alpha_{t-1i_{t-1}} \beta_{t-1},
\end{equation}
where $\alpha_{ij} \in C^*(u^{k_{ij}}p_1u^{k_{ij}}, ..., u^{k_{ij}}p_mu^{k_{ij}})^{\circ}$ with 
$0 \leq k_{ij} \leq n-1$, $k_{ij} \neq k_{i(j+1)}$ and $\beta_i \in C^*(e_{11}, ..., e_{nn})^{\circ}$. \\ 
\par

We consider the following two cases: \\ 

(a) We look at $\alpha_{ji} \alpha_{ji+1}$ with $\alpha_{jc}$ $\in$ $C^*(\{ u^{k_{c}}p_1u^{-k_{c}}, ... , 
u^{k_{c}}p_mu^{-k_{c}} \} )^{\circ}$ for 
$c=i, i+1$. 
We write $\alpha_{jc} = u^{k_{c}} \alpha'_{jc} u^{-k_{c}}$ with $\alpha'_{jc} 
\in C^*( \{ p_1, ... , p_m \} )^{\circ}$ for 
$c = i, i+1$. So $\alpha_{ji} \alpha_{ji+1} =$ \\ 
$u^{k_i} \alpha'_{ji} u^{k_{i+1} - k_i} \alpha'_{ji+1} u^{-k_{i+1}}$. Here $\alpha'_{ji}$ and $\alpha'_{ji+1}$ are free from 
$u^{k_{i+1} - k_i}$ in $(A,\tau)$ (Notice that we have $k_{i+1} - k_i \neq 0$). \\ 

(b) We look at $\alpha_{ji_j} \beta_j \alpha_{(j+1) 1}$ with $\beta \in C^*( \{e_{11} , ... , e_{nn} \}
)^{\circ},$ \\ 
$\alpha_{(j+1)1} \in C^*( \{ u^{k_{j+1}}p_1u^{-k_{j+1}} , ... , u^{k_{j+1}} p_m u^ {-k _{j+1}} \} )^{\circ},$
\\ 
$\alpha_{ji_j} \in C^*( \{u^{k_j} p_1 u^{-k_j} , ... , u^{k_j} p_m
u^{-k_j} \} )^{\circ}$. Now we write $\alpha_{ji_j} = u^{k_j} \alpha'_{ji_j}
 u^{-k_j}$ and $\alpha_{(j+1)1} = u^{k_{j+1}} \alpha'_{(j+1)1} u^{-k_{j+1}}$ with 
 $\alpha'_{ji_j} , \alpha'_{(j+1)1} \in 
 C^*( \{ p_1 , ..., p_m \} )^{\circ}$. We see that $\alpha_{ji_j} \beta_j \alpha_{(j+1)1} =$ $u^{k_j}
\alpha'_{ji_j} u^{-k_j} \beta_j u^{k_{j+1}} \alpha'_{(j+1)1} u^{-k_{j+1}}$. If $k_j = k_{j+1}$ then $\tau(u^{-k_j} \beta_j u^{k_{j+1}})
= \tau(u^{k_{j+1}} u^{-k_j} \beta_j) = \tau(\beta_j) = 0$ since $\tau$ is a trace. If $k_j \neq k_{j+1}$ then $\tau(u^{-k_j} \beta_j
u^{k_{j+1}}) = \tau(u^{k_{j+1}} u^{-k_j} \beta_j)$ and $u^{k_{j+1} - k_j} \beta_j \in \mathbb{M}_n$ is a linear combination of off-diagonal
elements, so $\tau(u^{k_{j+1}} u^{-k_j} \beta_j) = 0$ also. Notice that $\alpha'_{ji_j}$ and $\alpha'_{(j+1)1}$ are free from 
$u^{-k_j} \beta_j u^{k_{j+1}}$ in $(A,\tau)$. \\ 
 
 Now we expand all the letters in the word $\omega$ according to the cases (a) and (b). We
see that we obtain a word, consisting of letters of zero trace, such that every two consequitive letters come either from $C^*( \{p_1, ..., p_m \} )$  or from $\mathbb{M}_n$. So $\tau(\omega) = 0$. It only remains to look at the case of the word 
$\omega u^r$ 
which is the word $\omega$, but ending in $u^r$. There are two principally different cases for $\omega u^r$ 
from the all four possible choices for $\omega$: \\
 
In cases (1) and (2) $\alpha_{ti_t} = u^k \alpha'_{ti_t} u^{-k}$ for some $0 \leq k \leq n-1$  with 
$\alpha'_{ti_t} \in C^*( \{ p_1 , ..., p_m \} )^{\circ}$. So the word will end in $u^k \alpha'_{ti_t} u^{r-k}$. 
If $r = k$ then $\alpha'_{ti_t}$ will be the last letter with trace zero and everything else will be the 
same as for $\omega$, so the whole word will have trace $0$. If $k \neq r$ then $\tau(u^{r-k}) = 0$ and $u^{r-k}$ 
is free from $\alpha'_{ti_t}$ so the word in this case will be of zero trace too. \\

In cases (3) and (4) if $\beta_{t-1} u^{r}$ is the whole word then $\beta_{t-1} u^{r}$ is a linear 
combination of off-diagonal elements of $\mathbb{M}_n$, and so its trace is $0$. If not then 
$\alpha_{(t-1)i_{t-1}} = u^k \alpha'_{(t-1)i_{t-1}} u^{-k}$ with $\alpha'_{(t-1)i_{t-1}} \in 
C^*( \{ p_1 , ... , p_m \} )^{\circ}$. So the word ends in \\
$ u^k \alpha'_{(t-1)i_{t-1}} u^{-k} \beta_{t-1} 
u^{r} $. Similarly as above we see that $\tau(u^{-k} \beta_{t-1} u^{r}) = 0$ for all values of $k$ and $r$. 
The rest of the word we treat as above and conclude that it's of zero trace in this case too. \\

So in all cases $\tau( \omega u^r) = 0$ just what we had to show. \\ 
\par
Case (ii): \\
\par
As in case (i) 
$$\omega \in  \Lambda^{\circ}(C^*(p_1,..., p_m)^{\circ},..., C^*(u^{l-1}p_1u^{1-l},..., 
u^{l-1}p_mu^{1-l})^{\circ}, C^*( e_{11},..., e_{nn}, u^l,..., u^{n-l})^{\circ})$$ 
is of one of the following types: \\

\begin{equation}
\omega =  \alpha_{11} \alpha_{12} \cdots \alpha_{1i_1} \beta_1 \alpha_{21} \cdots \alpha_{2i_2} \beta_2
\alpha_{31} \cdots \alpha_{t-1i_{t-1}} \beta_{t-1} \alpha_{t1} \cdots \alpha_{ti_t},
\end{equation}

\begin{equation}
\omega =  \beta_1 \alpha_{21} \cdots \alpha_{2i_2} \beta_2
\alpha_{31} \cdots \alpha_{t-1i_{t-1}} \beta_{t-1} \alpha_{t1} \cdots \alpha_{ti_t},
\end{equation}

\begin{equation}
\omega =  \beta_1 \alpha_{21} \cdots \alpha_{2i_2} \beta_2
\alpha_{31} \cdots \alpha_{t-1i_{t-1}} \beta_{t-1},
\end{equation}

\begin{equation}
\omega =  \alpha_{11} \alpha_{12} \cdots \alpha_{1i_1} \beta_1 \alpha_{21} \cdots \alpha_{2i_2} \beta_2
\alpha_{31} \cdots \alpha_{t-1i_{t-1}} \beta_{t-1},
\end{equation}
where $\alpha_{ij} \in C^*(u^{k_{ij}}p_1u^{k_{ij}}, ..., u^{k_{ij}}p_mu^{k_{ij}})^{\circ}$ with 
$0 \leq k_{ij} \leq l-1$ and $k_{ij} \neq k_{(i+1)j}$ and $\beta_i \in C^*(e_{11}, ..., e_{nn}, u^l, 
u^{2l}, ..., u^{n-l})^{\circ}$. \\ 
\par  
Similarly as case (i) we consider two cases: \\ 

(a)   We look at $\alpha_{ji} \alpha_{ji+1}$ with $\alpha_{jc}$ $\in$ $C^*(\{ u^{k_{c}}p_1u^{-k_{c}}, ... ,
u^{k_{c}}p_mu^{-k_{c}} \} )$, and $0 \leq k_c \leq l-1$ for $c=i, i+1$. We write 
$\alpha_{jc} = u^{k_c} \alpha'_{jc} u^{-k_c}$ with $\alpha'_{jc} \in 
C^*( \{ p_1, ... , p_m \} )^{\circ}$ for $c = i, i+1$. It follows $\alpha_{ji} \alpha_{ji+1} =$ 
$u^{k_i} \alpha'_{ji} u^{k_{i+1} - k_i} \alpha'_{ji+1} u^{-k_{i+1}}$. Here $\alpha'_{ji}$ and $\alpha'_{ji+1}$ are free from $u^{k_{i+1} - k_i}$ in $(A,\tau)$ (and again $k_{i+1} - k_i \neq 0$). \\ 

(b)  We look at $\alpha_{ji_j} \beta_j \alpha_{(j+1) 1}$ with $\beta_j \in C^*( \{e_{11} , ... , e_{nn} \} , 
\{ u^l, u^{2l}, ..., u^{n-l} \} )^{\circ},$ \\ 
$\alpha_{(j+1)1} \in C^*( \{ u^{k_{j+1}}p_1u^{-k_{j+1}} , ... , u^{k_{j+1}} p_m u^ {-k _{j+1}} \} )^{\circ},$
\\ 
$\alpha_{ji_j} \in C^*( \{u^{k_j} p_1 u^{-k_j} , ... , u^{k_j} p_m u^{-k_j} \} )^{\circ}$, where in this case 
$k_j, k_{j+1} \in \{ 0, ..., l-1 \}$. Again we write $\alpha_{ji_j} = u^{k_j} \alpha'_{ji_j} u^{-k_j}$ and 
$\alpha_{(j+1)1} = u^{k_{j+1}} \alpha'_{(j+1)1} u^{-k_{j+1}}$ 
with $\alpha'_{ji_j} , \alpha'_{(j+1)1} \in C^*( \{ p_1 , ... , p_m \} )^{\circ},$. 
We have $\alpha_{ji_j} \beta_j \alpha_{(j+1)1} =$ $u^{k_j} \alpha'_{ji_j} u^{-k_j} \beta_j u^{k_{j+1}} \alpha'_{(j+1)1} u^{-k_{j+1}}$. \\ 
We only need to show that $\tau(u^{-k_j} \beta_j u^{k_{j+1}}) = 0$. $\tau(u^{-k_j} \beta_j u^{k_{j+1}}) = \tau(u^{k_{j+1}} u^{-k_j} 
\beta_j) =  \tau(u^{k_{j+1} - k_j} \beta_j)$. The case $ k_{j+1} = k_j$ is clear. Notice that if 
$ k_{j+1} \neq k_j$ then $0 < k_{j+1} - k_j \leq l-1$. Is it clear that 
$u^{k_{j+1} - k_j} \cdot \Span ( \{ e_{11}, ..., e_{nn} \}) \subset \mathbb{M}_n$ consists of liner 
combination of off-diagonal elements. The same is clear for $u^{k_{j+1} - k_j} \cdot 
\Span( \{ u^l, u^{2l} , ..., u^{n-l} \} ) \subset \mathbb{M}_n $. It's not difficult to see then that 
$$u^{k_{j+1} - k_j} \cdot \Alg ( \{ e_{11}, ..., e_{nn} \}, \{ u^l, u^{2l}, ..., u^{n-l} \} )$$ 
will consist of linear span of the union of the off-diagonal entries among 
$\{ e_{ij} | 1 \leq i,j \leq n \}$ present in $u^{k_{j+1} - k_j} \cdot \Span( \{e_{11}, ..., e_{nn} \})$ and 
the ones present in \\ 
$u^{k_{j+1} - k_j} \cdot \Span( \{ u^l, u^{2l}, ..., u^{n-l} \} )$. 
This shows that $u^{k_{j+1} - k_j} \beta_j$ will be also a linear span of off-diagonal entries in 
$\mathbb{M}_n$ and will have trace $0$. 
So $\tau(u^{-k_j} \beta_j u^{k_{j+1}}) = 0$.  In this case also $\alpha'_{ji_j}$ and $\alpha'_{(j+1)1}$ are 
free from $u^{-k_j} \beta_j u^{k_{j+1}}$ in $(A,\tau)$. \\ 

We expand all the letters of the word $\omega$ and see that it is of trace $0$ similarly as in case (i). 
For the word $\omega u^r$ with $0 \leq r \leq l-1$ we argue similarly as in case (i). 
Again there are two principally different cases: \\

In cases (5) and (6) $\alpha_{ti_t} = u^k \alpha'_{ti_t} u^{-k}$ for some $0 \leq k \leq l-1$  with 
$\alpha'_{ti_t} \in C^*( \{ p_1 , ..., p_m \} )^{\circ}$. So the word will end in $u^k \alpha'_{ti_t}
u^{r-k}$. If $r = k$ then $\alpha'_{ti_t}$ will be the last letter with trace zero and everything else will 
be the same as for $\omega$, so the whole word will have trace $0$. If $k \neq r$ then $\tau(u^{r-k}) = 0$ 
and $u^{r-k}$ is free from $\alpha'_{ti_t}$ so the word in this case will be of zero trace too.

In cases (7) and (8) $\beta_{t-1} u^r$ then this is a linear combination of off-diagonal elements as we showed 
in case (ii)-(b). If not we write $\alpha_{(t-1)i_{t-1}} = u^k \alpha'_{(t-1)i_{t-1}} u^{-k}$ with 
$0 \leq k \leq l-1$ and $\alpha'_{(t-1)i_{t-1}} \in 
C^*( \{ p_1, ..., p_m \} )^{\circ}$. So the word that we are looking at will end in 
$u^k \alpha'_{(t-1)i_{t-1}} u^{-k} \beta_{t-1} u^{r} $. Since $0 \leq k,r \leq l-1$ similarly as in case
(ii)-(b) we see that $\tau(u^{-k} \beta_{t-1} u^{r}) = 0$. We treat the remaining part of the word as above and conclude that in this case the word has trace 
$0$. \\ 
\par
So in all cases $\tau(\omega u^r) = 0$ just what we had to show. \\
\par
This proves the lemma.

\end{proof}

From properties (5) and (6) of the reduced free product it follows that $\tau$ is a faithful trace. From Lemma 1.4 it follows that $$B =
(\mathbb{C} \cdot e_{11} \oplus ... \oplus \mathbb{C} \cdot e_{nn}) * (\underset{k=0}{\overset{n-1}{*}} (\mathbb{C} \cdot u^k p_1 u^{-k} \oplus ... \oplus
\mathbb{C} \cdot u^k p_m u^{-k})),$$ 
$$\cong ( \underset{\frac{1}{n}}{\mathbb{C}} \oplus ... \oplus \underset{\frac{1}{n}}{\mathbb{C}} ) * 
(\underset{k=0}{\overset{n-1}{*}} (\underset{\alpha_1}{\mathbb{C}} \oplus ... \oplus
\underset{\alpha_m}{\mathbb{C}}))$$ and that $$E = 
C^*( \{ e_{11}, ..., e_{nn}, u_l, u^{2l}, ..., u^{n-l} \} ) * (\underset{k=0}{\overset{l-1}{*}} (\mathbb{C} \cdot u^k p_1 u^{-k} \oplus ... \oplus
\mathbb{C} \cdot u^k p_m u^{-k})),$$ 
$$\cong (\underset{\frac{l}{n}}{\mathbb{M}_{\frac{n}{l}}} \oplus ... \oplus 
\underset{\frac{l}{n}}{\mathbb{M}_{\frac{n}{l}}}) * (\underset{k=0}{\overset{l-1}{*}} (\underset{\alpha_1}{\mathbb{C}} \oplus ...
\oplus \underset{\alpha_m}{\mathbb{C}})).$$ 

\begin{cor}

If everything is as above:

(1)  For $b \in B$ and $0 < k \leq n-1$ we have $\tau(b u^k) = 0$, so also $\tau(u^k b) = 0$. 
\par
(2)  For $e \in E$ and $0 < k \leq l-1$ we have $\tau(e u^k) = 0$, so also $\tau(u^k e) = 0$.

\end{cor}

For $(B, \tau|_B)$ and $(E, \tau|_E)$ we have that $\mathfrak{H}_B \subset \mathfrak{H}_E \subset \mathfrak{H}_A$. If $a \in A$ we will 
denote by $\hat{a} \in \mathfrak{H}_A$ the vector in $\mathfrak{H}_A$, corresponding to $a$ by the GNS construction. We will show that
 
\begin{cor}

If everything is as above: \\ 

(1)  $u^{k_1} \mathfrak{H}_B  \bot u^{k_2} \mathfrak{H}_B$ for $k_1 \neq k_2$, $0 \leq k_1, k_2 \leq n-1$. 
\par
(2)  $u^{k_1} \mathfrak{H}_E  \bot u^{k_2} \mathfrak{H}_E$ for $k_1 \neq k_2$, $0 \leq k_1, k_2 \leq l-1$.

\end{cor}

\begin{proof} 

(1)  Take $ b_1, b_2 \in B $. We have $\langle u^{k_1} \hat{b_1} , u^{k_2} \hat{b_2} \rangle = \tau(u^{k_2} b_2 b_1^* u^{-k_1}) = 
\tau(b_2 b_1^* u^{k_2 - k_1}) = 0,$ by the above Corollary. \\ 
(2)  Similarly take $e_1, e_2 \in E$, so $\langle u^{k_1} \hat{e_1}, u^{k_2} \hat{e_2} \rangle = \tau(u^{k_2} e_2 e_1^* u^{-k_1}) 
= \tau(e_2 e_1^* u^{k_2 - k_1}) = 0,$ again by the above Corollary.

\end{proof}

Now $\mathfrak{H}_A$ can be written in the form $\mathfrak{H}_A = 
\underset{i=0}{\overset{n-1}{\bigoplus}} u^i \mathfrak{H}_B$ as a
Hilbert space because of the Corollary above. Denote by $P_i$ the projection $P_i : \mathfrak{H}_A \rightarrow \mathfrak{H}_A$ onto the
subspace $u^i \mathfrak{H}_B$. Now it's also true that $A = 
\underset{i=0}{\overset{n-1}{\bigoplus}} u^i B$ as a Banach
space. To see this we notice that $\Span\{u^iB, i=0, ...n-1\}$ is dense in $A$, also that $u^i B,\ 0 \leq i 
\leq n-1$ are closed in $A$. Now take a sequence $\{ \sum_{i=0}^{n-1} u^i b_{mi} \}_{m=1}^{\infty}$ converging to an element $a \in A$ 
($b_{mi} \in B$). 
Then for each $i$ we have $\{ P_j  \sum_{i=0}^{n-1} u^i b_{mi} P_0 \}_{m=1}^{\infty}  = \{ P_j u^j b_{mj} P_0 \}_{m=1}^{\infty}$ converges 
(to $P_j a P_0$), consequently the
sequence $\{ b_{mj} \}_{m=1}^{\infty}$ converges to an element $b_j$ in $B$ $\forall 0 \leq j \leq n-1$. So 
$a = \sum_{i=0}^{n-1} u^i b_i$.  Finally 
the fact that $u^{i_1} B \cap u^{i_2} B = 0$, for $i_1 \neq i_2$ follows easily from 
$u^{i_1} \mathfrak{H}_B \cap u^{i_2} \mathfrak{H}_B = 0$, for $i_1 \neq i_2$ and the fact that the trace $\tau$ is faithful. We also have
$A = \underset{i=0}{\overset{n-1}{\bigoplus}} B u^i$.
\par
Let $C$ is a $C^*$-algebra and $\Gamma$ is a discrete group with a given action 
$\alpha : \Gamma \rightarrow Aut(C)$ on $C$. By $C \rtimes \Gamma$ we will denote the reduced crossed 
product of $C$ by $\Gamma$. It will be clear what group action we take.
\par
Let's denote by $G$ the multiplicative group, generated by the automorphism $\Ad(u)$ of $B$. Then 
$G \cong \mathbb{Z}_n$ and by what we proved above 
$\mathfrak{H}_A \cong L^2(G,\mathfrak{H}_B)$. 

\begin{lemma}

$A \cong B \rtimes G$

\end{lemma}

\begin{proof}

We have to show that the action of $A$ on $\mathfrak{H}_A$ "agrees" with the crossed product action. Take $a=
\underset{k=0}{\overset{n-1}{\sum}} b_k u^k \in A$, $b_k \in B, k=0, 1, ..., n-1$ and take $\xi = \underset{k=0}{\overset{n-1}{\sum}} u^k
\hat{b'_k} \in \mathfrak{H}_A$, $b'_k \in B, k=0, 1, ..., n-1$. Then $$a(\xi) = \underset{k=0}{\overset{n-1}{\sum}} 
\underset{m=0}{\overset{n-1}{\sum}} b_k u^k u^m \hat{b'_m} = \underset{k=0}{\overset{n-1}{\sum}} \underset{m=0}{\overset{n-1}{\sum}}
u^{k+m} . (u^{-k-m} b_k u^{k+m} ) \hat{b'_m},$$ 
$$= \underset{s=0}{\overset{n-1}{\sum}} \underset{k=0}{\overset{n-1}{\sum}} (u^s . \Ad(u^{-s})(b_k) ) (\widehat{b'_{s-k(mod\ n)}}).$$ This shows that the action of $A$ on $\mathfrak{H}_A$ is the crossed product action.

\end{proof}

To study simplicity in this situation, we can invoke Theorem 4.2 from \cite{O75} and Theorem 6.5 from \cite{OP78}, or with the same success, use the following result from \cite{K81}:

\begin{thm}[\cite{K81}]

Let $\Gamma$ be a discrete group of automorphisms of $C^*$-algebra $\mathfrak{B}$. If $\mathfrak{B}$ is simple and if each $\gamma$ is 
outer for the multiplier algebra $M(\mathfrak{B})$ of $\mathfrak{B}$, $\forall \gamma \in \Gamma \backslash \{ 1 \} $, then the
reduced crossed product of $\mathfrak{B}$ by $\Gamma$, $\mathfrak{B} \rtimes \Gamma$, is simple.

\end{thm}

An automorphism $\omega$ of a $C^*$-algebra $\mathfrak{B}$ , contained in a $C^*$-algebra $\mathfrak{A}$ is 
outer for $\mathfrak{A}$, if there doesn't exist a unitary $w \in \mathfrak{A}$ with the property $\omega = \Ad(w)$. 
\par
A representation $\pi$ of a $C^*$-algebra $\mathfrak{A}$ on a Hilbert space $\mathfrak{H}$ is called 
non-degenerate if there doesn't exist a vector $\xi \in
\mathfrak{H}$, $\xi \neq 0$, such that $\pi(\mathfrak{A}) \xi = 0$.
\par
The idealizer of a $C^*$-algebra $\mathfrak{A}$ in a $C^*$-algebra $\mathfrak{B}$ ($\mathfrak{A} \subset \mathfrak{B}$) is the largest $C^*$-subalgebra of $\mathfrak{B}$ in which $\mathfrak{A}$ is an ideal. \\ 
We will not give a definition of multiplier algebra of a $C^*$-algebra. Instead we will give the following property from \cite{APT73}, which we will use (see \cite{APT73} for more details on multiplier algebras):

\begin{prop}[\cite{APT73}]

Each nondegenerate faithful representation $\pi$ of a $C^*$-algebra $\mathfrak{A}$ extends uniquely to a faithful representation of
$M(\mathfrak{A})$, and $\pi(M(\mathfrak{A}))$ is the idealizer of $\pi(\mathfrak{A})$ in its weak closure.

\end{prop}

Suppose that we have a faithful representation $\pi$of a $C^*$ algebra $\mathfrak{A}$ on a Hilbert space 
$\mathfrak{H}$. If confusion is impossible we will denote by $\bar{\mathfrak{A}}$ (in $\mathfrak{H}$) the 
weak closure of $\pi(\mathfrak{A})$ in $\mathbb{B}(\mathfrak{H})$.
\par
To study uniqueness of trace we invoke a theorem of B$\acute{e}$dos from \cite{B93}. 
\par
Let $\mathfrak{A}$ be a simple, unital $C^*$-algebra with a unique trace $\varphi$ and let 
$(\pi_{\mathfrak{A}}, \mathfrak{H}_{\mathfrak{A}}, 
\widehat{1_{\mathfrak{A}}})$ denote the GNS-triple associated to $\varphi$. The trace $\varphi$ is faithful by the simplicity of $\mathfrak{A}$ and 
$\mathfrak{A}$ is isomorphic to $\pi_{\mathfrak{A}}(\mathfrak{A})$. Let
$\alpha \in Aut(\mathfrak{A})$. The trace $\varphi$ is $\alpha$-invariant by the uniqueness of $\varphi$. Then $\alpha$ is implemented on 
$\mathfrak{H}_{\mathfrak{A}}$ by the unitary operator $U_{\alpha}$ given by 
$U_{\alpha}(\hat{a}) = \alpha(a) \cdot \widehat{1_{\mathfrak{A}}}$, 
$a \in \mathfrak{A}$.  Then we denote the extension of $\alpha$ to the weak closure 
$\bar{\mathfrak{A}}$ (in $\mathfrak{H}_{\mathfrak{A}}$) of $\pi_{\mathfrak{A}}(\mathfrak{A})$ on $\mathbb{B}(\mathfrak{H}_{\mathfrak{A}})$ 
by $\tilde{\alpha} \overset{def}{=} \Ad(U_{\alpha})$. We will say that $\alpha$ is $\varphi$-outer if $\tilde{\alpha}$ is outer for $\bar{\mathfrak{A}}$. 

\begin{thm}[\cite{B93}]

Suppose $\mathfrak{A}$ is a simple unital $C^*$-algebra with a unique trace $\varphi$ and that $\Gamma$ is a discrete group with a 
representation $\alpha : \Gamma \rightarrow Aut(\mathfrak{A})$, such that $\alpha_{\gamma}$ is $\varphi$-outer $\forall \gamma \in \Gamma 
\backslash \{ 1 \}$. Then the reduced crossed product $\mathfrak{A} \rtimes \Gamma$ is simple with a unique trace $\tau$ given by 
$\tau = \varphi \circ E$, where $E$ is the canonical conditional expectation from $\mathfrak{A} \rtimes \Gamma$ onto $\mathfrak{A}$.

\end{thm}

Let's now return to the $C^*$-algebra $(A,\tau) =  ( \underset{\alpha_1}{\overset{p_1}{\mathbb{C}}} \oplus ... \oplus 
\underset{\alpha_m}{\overset{p_m}{\mathbb{C}}})*(\mathbb{M}_n, tr_n)$, with $\alpha_1 \leq \alpha_2 \leq ... \leq \alpha_m$. If 
$B \subset E \subset A$ are as in the beginning of this section, then the representations of $B$, $E$ and $A$ on $\mathfrak{H}_A$ are all nondegenerate. Also we have the following:

\begin{lemma}

The weak closure of $B$ in $\mathbb{B}(\mathfrak{H}_B)$ and the one in $\mathbb{B}(\mathfrak{H}_A)$ are the 
same (or $\bar{B}$ (in $\mathfrak{H}_B$) $\cong$ 
$\bar{B}$ (in $\mathfrak{H}_A$)). Analoguously, $\bar{E}$ (in $\mathfrak{H}_E$) $\cong$ $\bar{E}$ (in $\mathfrak{H}_A$).

\end{lemma}

\begin{proof}

For $b \in B \subset A$ we have $b(u^t h) = u^t (\Ad(u^{-t} b))(h)$ for $h \in \mathfrak{H}_B$ and 
$0 \leq t \leq n-1$. Taking a weak limit
in $\mathbb{B}(\mathfrak{H}_B)$ we obtain the same equation $\forall \bar{b} \in \bar{B}$ 
(in $\mathfrak{H}_B$): $\bar{b}(u^th) = 
u^t(\Ad(u^{-t})(\bar{b}))(h)$, which shows, of course, that $\bar{b}$ has a unique extension to 
$\mathbb{B}(\mathfrak{H}_A)$. Conversely if $\tilde{b} \in \bar{B}$ (in
$\mathfrak{H}_A$), then since $\mathfrak{H}_B$ is invariant for $B$ it will be invariant for $\tilde{b}$ also. So the restriction of
$\tilde{b}$ to $\mathfrak{H}_B$ is the element we are looking for. 
\par
Analoguously if $e \in E$ and if $h_0 + u^l h_1 + ... + u^{n-l} h_{\frac{n}{l}-1} \in \mathfrak{H}_E$, then for $0 \leq t \leq l-1$ we have
$e(u^t(h_0 + u^l h_1 + ... + u^{n-l} h_{\frac{n}{l}-1})) = u^t(\Ad(u^{-t})(e))(h_0 + u^l h_1 + ... + u^{n-l} h_{\frac{n}{l}-1})$. And again
for an element $\bar{e} \in \bar{E}$ (in $\mathfrak{H}_E$) we see that $\bar{e}$ has a unique extension to an element of $\bar{E}$ (in
$\mathfrak{H}_A$). Conversely an element $\tilde{e} \in \bar{E}$ (in $\mathfrak{H}_A$) has $\mathfrak{H}_E$ as an invariant subspace, so we 
can restrict it to $\mathfrak{H}_E$ to obtain an element in $\bar{E}$ (in $\mathfrak{H}_E$).

\end{proof}

We will state the following theorem from \cite{D99}, which we will frequently use:

\begin{thm}[\cite{D99}]

Let $\mathfrak{A}$ and $\mathfrak{B}$ be unital $C^*$-algebras with traces $\tau_{\mathfrak{A}}$ and
$\tau_{\mathfrak{B}}$ respectively, whose GNS representations are faithful. Let 

\begin{center}

$(\mathfrak{C}, \tau) = (\mathfrak{A}, \tau_{\mathfrak{A}}) * (\mathfrak{B}, \tau_{\mathfrak{B}})$.

\end{center}

Suppose that $\mathfrak{B} \neq \mathbb{C}$ and that $\mathfrak{A}$ has a unital, diffuse abelian 
$C^*$-subalgebra $\mathfrak{D}$ ($1_{\mathfrak{A}} \in \mathfrak{D} \subseteq \mathfrak{A}$). 
Then $\mathfrak{C}$ is simple with a unique trace $\tau$.

\end{thm}

Using repeatedly Theorem 2.4 we see that $$B =
(\mathbb{C} \cdot e_{11} \oplus ... \oplus \mathbb{C} \cdot e_{nn}) * (\underset{k=0}{\overset{n-1}{*}} (\mathbb{C} \cdot u^k p_1 u^{-k} \oplus ... \oplus
\mathbb{C} \cdot u^k p_m u^{-k})),$$ 
$$\cong (U \oplus \underset{max \{ n\alpha_m - n + 1,\ 0 \} }{\overset{\tilde{p}}{\mathbb{C}}}) *
(\underset{\frac{1}{n}}{\overset{e_{11}}{\mathbb{C}}} \oplus ... \oplus \underset{\frac{1}{n}}{\overset{e_{nn}}{\mathbb{C}}}),$$ where $U$ 
has a unital, diffuse abelian $C^*$-subalgebra, and where $\tilde{p} = \underset{i=0}{\overset{n-1}{\wedge}} u^i p_m u^{-i}$. 
\par
We will consider the following 3 cases, for $\alpha_1 \leq \alpha_2 \leq ... \leq \alpha_m$: \\ 
\par
(I) $\alpha_m < 1-\frac{1}{n^2}$. 
\par
(II) $\alpha_m = 1-\frac{1}{n^2}$. 
\par
(III) $\alpha_m > 1-\frac{1}{n^2}$. \\ 
\par
We will organize those cases in few lemmas:
\par
\par

(I)

\begin{lemma}

If $A$ is as above, then for $\alpha_m < 1-\frac{1}{n^2}$ we have that $A$ is simple with a unique trace.

\end{lemma} 

\begin{proof}

We consider: \\ 
(1) $\alpha_m \leq 1-\frac{1}{n}$. \\
Then $B \cong U * (\underset{\frac{1}{n}}{\overset{e_{11}}{\mathbb{C}}} \oplus ... \oplus 
\underset{\frac{1}{n}}{\overset{e_{nn}}{\mathbb{C}}})$ with $U$ containing a unital, diffuse abelian $C^*$-subalgebra (from Theorem 2.4). 
From the Theorem 3.9 we see that $B$ is simple with a unique trace. \\ 
(2) $1-\frac{1}{n} < \alpha_m < 1-\frac{1}{n^2}$. \\ 
Then $B \cong (U \oplus \underset{n\alpha_m - n + 1}{\overset{\tilde{p}}{\mathbb{C}}}) *
(\underset{\frac{1}{n}}{\overset{e_{11}}{\mathbb{C}}} \oplus ... \oplus \underset{\frac{1}{n}}{\overset{e_{nn}}{\mathbb{C}}})$ with $U$
having a unital, diffuse abelian $C^*$-subalgebra. Using Theorem 2.4 one more time we see that $B$ is simple with a unique trace in this case also.
\par
We know that $A = B \rtimes G$, where $G = \langle \Ad(u) \rangle \cong \mathbb{Z}_n$. Since $B$ is unital then the multiplier algebra $M(B)$
coinsides with $B$. We note also that since $\bar{B}$ (in $\mathfrak{H}_B$)
is isomorphic to $\bar{B}$ (in $\mathfrak{H}_A$) to prove that some element of $Aut(B)$ is $\tau_B$-outer it's enough to prove that this
automorphism is outer for $\bar{B}$ (in $\mathfrak{H}_A$) (and it will be outer for $M(B) = B$ also). Making these observations and using 
Theorem 3.5 and Theorem 3.7 we see that if we prove that $\Ad(u^i)$ is outer for $\bar{B}$ (in $\mathfrak{H}_A$), $\forall 0 < i \leq n-1$, 
then it will follow that $A$ is simple with a unique trace. We will show that $\Ad(u^i)$ is outer for $\bar{B}$ (in $\mathfrak{H}_A$) 
(we will write just $\bar{*}$ for $\bar{*}$ (in $\mathfrak{H}_A$) and omit writting $\mathfrak{H}_A$ - all the
closures will be in $\mathbb{B}(\mathfrak{H}_{\mathfrak{A}})$) for 
the case  $\alpha_m \leq 1-\frac{1}{n^2}$. 
\par
Fix $0 < k \leq n-1$. Since $u^k \mathfrak{H}_B \perp \mathfrak{H}_B$ it follows that $u^k \notin \bar{B}$ (in $\mathfrak{H}_A$). Suppose
$\exists w \in \bar{B}$, such that $\Ad(u^k) = \Ad(w)$ on $\bar{B}$. Then $u^k w u^{-k} = w w w^* = w$ and $u^k w^* u^{-k} = w w^* w^* = w^*$ 
and this implies that $u^k$, $u^{-k}$, $w$ and $w^*$ commute, so it follows $u^k w^*$ commutes with $\overline{C^*(B, u^k)}$, so it belongs
to its center. If $k \nmid n$ then $\overline{C^*(B, u^k)} = \bar{A}$ and by Theorem 2.5 $\bar{A}$ (in $\mathfrak{H}_A$)is a factor, so $u^k w^*$ is a 
multiple of $1_A$, which contradicts the fact $u^k \notin \bar{B}$. If $k=l \mid n$, then $\overline{C^*(B, u^k)} = \bar{E}$ and $\bar{E}$ 
(in $\mathfrak{H}_A$) $\cong$ $\bar{E}$ (in $\mathfrak{H}_E$) is a factor too (by Theorem 2.5), so this implies again that $u^k w^*$ is a 
multiple of $1_A = 1_E$, so this is a contradiction again and this proves that $\Ad(u^k)$ are outer for $\bar{B}$, $\forall 0 < k \leq n-1$.
This concludes the proof.

\end{proof}

(III) 

\begin{lemma}

If $A$ is as above, then for $\alpha_m > 1-\frac{1}{n^2}$ we have $A = A_0 \oplus \underset{n^2 \alpha_m - n^2 + 1}{\mathbb{M}_n}$, where $A_0$ is simple with a unique trace.

\end{lemma}

\begin{proof}

In this case $B \cong (U \oplus \underset{n\alpha_m - n + 1}{\overset{\tilde{p}}{\mathbb{C}}}) * 
(\underset{\frac{1}{n}}{\overset{e_{11}}{\mathbb{C}}} \oplus ... \oplus \underset{\frac{1}{n}}{\overset{e_{nn}}{\mathbb{C}}})$, where $U$
has a unital, diffuse abelian $C^*$-subalgebra. Form Theorem 2.4 we see that $B \cong \overset{\tilde{p}_0}{B_0} \oplus \underset{n
\alpha_m - n + \frac{1}{n}}{\overset{e_{11} \wedge \tilde{p}}{\mathbb{C}}} \oplus ... \oplus \underset{n
\alpha_m - n + \frac{1}{n}}{\overset{e_{nn} \wedge \tilde{p}}{\mathbb{C}}}$ with $\tilde{p}_0 = 1- e_{11} \wedge 
\tilde{p} - ... - e_{nn} \wedge \tilde{p}$, and $B_0$ being a unital, simple and having a
unique trace. It's easy to see that $\Ad(u)$ permutes $\{ e_{ii} | 1 \leq i \leq n \}$ and that $\Ad(u)$ permutes 
$\{ u^i p_j u^{-i} | 0 \leq i \leq n-1 \}$ for each $1 \leq j \leq m$. But since $\tilde{p} = 
\underset{i=0}{\overset{n-1}{\wedge}} u^i p_m u^{-i}$ we see that $\Ad(u)(\tilde{p}) = \tilde{p}$. 
This shows that $\Ad(u)$ permutes 
$\{ e_{ii} \wedge \tilde{p} | 1 \leq i \leq n \}$. This shows that $\Ad(\tilde{p}_0 u)$ is an automorphism 
of $B_0$ and that 
$\Ad((1-\tilde{p}_0) u)$ is an automorphism of $\overset{e_{11} \wedge \tilde{p}}{\mathbb{C}} \oplus ... \oplus \overset{e_{nn} \wedge 
\tilde{p}}{\mathbb{C}}$. If we denote $G_1 = \langle \Ad(\tilde{p}_0 u) \rangle$ and $G_2 = \langle \Ad((1-\tilde{p}_0) u) \rangle$, 
then we have 
$A = B_0 \rtimes G_1 \oplus (\overset{e_{11} \wedge \tilde{p}}{\mathbb{C}} \oplus ... \oplus \overset{e_{nn} \wedge \tilde{p}}{\mathbb{C}})
\rtimes G_2$. Now it's easy to see that $(\overset{e_{11} \wedge \tilde{p}}
{\mathbb{C}} \oplus ... \oplus \overset{e_{nn} \wedge \tilde{p}}{\mathbb{C}}) \rtimes G_2 = C^*(\{ e_{11} \wedge \tilde{p}, ..., e_{nn} 
\wedge \tilde{p} \}, (1-\tilde{p}_0) u) = (1-\tilde{p}_0).C^*( \{ e_{11}, ..., e_{nn} \}, u) \cong \mathbb{M}_n$ (because $\tilde{p}_0$ is
a central projection). To study 
$A_0 \overset{def}{=} B_0 \rtimes G_1$ we have to consider the automorphisms $\Ad(\tilde{p}_0 u)$. From Lemma 3.8 we see that 
$$\overline{B_0 \oplus \overset{e_{11} 
\wedge\tilde{p}}{\mathbb{C}} \oplus ... \oplus \overset{e_{nn} \wedge \tilde{p}}{\mathbb{C}}}\ (in\ \mathfrak{H}_B) \cong \overline{B_0 
\oplus \overset{e_{11} \wedge\tilde{p}}{\mathbb{C}} \oplus ... \oplus \overset{e_{nn} \wedge \tilde{p}}{\mathbb{C}}}\ (in\ \mathfrak{H}_A).$$ 
This implies $\bar{B}_0$ (in $\mathfrak{H}_{B_0}$) $\cong$ $\bar{B}_0$ (in $\mathfrak{H}_{A_0}$). This is because 
$\mathfrak{H}_{A_0} = \tilde{p}_0 \mathfrak{H}_A$ and $\mathfrak{H}_{B_0} = \tilde{p}_0 \mathfrak{H}_B$ (which is clear, since
$\mathfrak{H}_{A_0}$ and $\mathfrak{H}_{B_0}$ are direct summands in $\mathfrak{H}_A$ and $\mathfrak{H}_B$ respectivelly). For some 
$l | n$ if we denote $E_0 \overset{def}{=} \tilde{p}_0 E$ then by the same reasoning as above $$E = E_0 \oplus (1-\tilde{p}_0).
C^*(\{ e_{11}, ..., e_{nn} \}, u^l) \cong E_0 \oplus (\underbrace{\mathbb{M}_{\frac{n}{l}} \oplus ... \oplus \mathbb{M}_{\frac{n}{l}}}_{l-times}).$$ 
So we similarly have $\bar{E_0}$ (in $\mathfrak{H}_{E_0}$) $\cong$ $\bar{E_0}$ (in $\mathfrak{H}_{A_0}$). We use Theorem 2.5 and see that 
$\bar{A} \cong L(F_t) \oplus \mathbb{M}_n$ and that $$\bar{E} \cong L(F_{t'}) \oplus (\underbrace{\mathbb{M}_{\frac{n}{l}} \oplus ... \oplus 
\mathbb{M}_{\frac{n}{l}}}_{l-times}),$$ for some $1 < t, t' < \infty$. This shows that $\bar{A_0}$ and $\bar{E_0}$ are both factors. 
Now for $\Ad(\tilde{p_0} u^k)$, $1 \leq k \leq n-1$ we can make the same reasoning as in the case (I) to 
show that $\Ad(\tilde{p_0} u^k)$ are all outer for 
$\bar{B_0}$, $\forall 1 \leq k \leq n-1$. Now we use Theorem 3.5 and Theorem 3.7 to finish the proof. Notice that the trace of the support
projection of $\mathbb{M}_n$, $e_{11} \wedge \tilde{p} + ... + e_{nn} \wedge \tilde{p}$, is $n^2 \alpha_m - n^2 + 1$.

\end{proof}

(II) \\ 
\par
We already proved that $\Ad(u^k)$ are outer for $\bar{B}$, $\forall 1 \leq k \leq n-1$. Using Theorem 2.4 we see 
$B \cong (U \oplus \underset{1-\frac{1}{n}}{\overset{\tilde{p}}{\mathbb{C}}}) * (\underset{\frac{1}{n}}{\overset{e_{11}}{\mathbb{C}}} 
\oplus ... \oplus \underset{\frac{1}{n}}{\overset{e_{nn}}{\mathbb{C}}})$ with $U$ having a unital, diffuse abelian $C^*$-subalgebra. There 
are $*$-homomorphisms $\pi_i : B \rightarrow \mathbb{C}$, $1 \leq i \leq n$ with $\pi_i(\tilde{p}) = \pi_i(e_{ii}) = 1$, 
and such that $B_0 \overset{def}{=} \underset{i=0}{\overset{n-1}{\bigcap}} \ker(\pi_i)$ is simple with a unique trace. Now if 
$1 \leq k \leq n-1$, then $B_0 \bigcap \Ad(u^k)(B_0) = $ either $0$ or $B_0$, because $B_0$ and $\Ad(u^k)(B_0)$ are simple ideals in $B$. 
The first possibility is actually 
impossible, because of dimension reasons, so this shows that $B_0$ is invariant for $\Ad(u^k)$, $1 \leq k \leq n-1$. In other words 
$\Ad(u^k) \in Aut(B_0)$. Similarly as in Lemma 3.4 it can be shown that 
$$A_0 \overset{def}{=} C^*(B_0 \oplus B_0 u \oplus ... \oplus B_0
u^{n-1}) \cong B_0 \rtimes \{ \Ad(u^k) | 0 \leq k \leq n-1 \} \subset A.$$ 

\begin{lemma}

We have a short split-exact sequence: 

\begin{center}

$0 \hookrightarrow A_0 \rightarrow A \overset{\curvearrowleft}{\rightarrow} \mathbb{M}_n \rightarrow 0$.

\end{center}

\end{lemma}

\begin{proof}

It's clear that we have the short exact sequence 

\begin{equation*}
0 \rightarrow B_0 \hookrightarrow B \overset{\pi}{\longrightarrow} \underbrace{\mathbb{C} 
\oplus ... \oplus \mathbb{C}}_{n-times} \rightarrow 0, 
\end{equation*}
where $\pi \overset{def}{=} (\pi_1, ..., \pi_n)$. We think $\pi$ to be a map from $B$ to $diag(\mathbb{M}_n)$, defined by $$\pi(b) = 
\begin{pmatrix} \pi_1(b) & 0 & ... & 0  \\ 0 & \pi_2(b) & ... & 0 \\   . & . & . & . \\ 0 & 0 & ... & \pi_n(b) \end{pmatrix} .$$ 
Now since 
$\pi_i(\tilde{p}) = \pi_i(e_{ii}) = 1$ and $\Ad(u)(e_{11}) = u e_{11} u^* = e_{nn}$ and for 
$2 \leq i \leq n$, $\Ad(u)(e_{ii}) = u e_{ii} 
u^* = e_{(i-1) (i-1)}$, then $\pi_i \circ \Ad(u)(e_{(i+1) (i+1)}) = \pi_i \circ \Ad(u)(\tilde{p}) = 1$ for $1 \leq i \leq n-1$ and $\pi_n 
\circ \Ad(u)(e_{1 1}) = \pi_n \circ \Ad(u)(\tilde{p}) = 1$. So since two $*$-homomorphism of a $C^*$-algebra, which coinside on a set 
of generators of the $C^*$-algebra, are identical, we have $\pi_i \circ \Ad(u) = \pi_{i+1}$ for 
$1 \leq i \leq n-1$ and $\pi_n \circ \Ad(u) = \pi_1$. Define $\tilde{\pi} : A \rightarrow \mathbb{M}_n$ by 
$\underset{k=0}{\overset{n-1}{\sum}} b_ku^k \mapsto 
\underset{k=0}{\overset{n-1}{\sum}} \pi(b_k) W^k$ (with $b_k \in B$), where $W \in \mathbb{M}_n$ is represented by the matrix, which 
represent $u \in \mathbb{M}_n \subset A$, namely $$W \overset{def}{=} \begin{pmatrix} 0 & 1 & ... & 0 \\ . & . & . & . \\ 0 & 0 & ... & 1 
\\  1 & 0 & ... & 0 \end{pmatrix} .$$
We will show that if $b \in B$ and $0 \leq k \leq n-1$, then $\pi(u^k b u^{-k}) = W^k \pi(b) W^{-k}$. For this it's enough to show that $\pi(u b u^{-1}) = W \pi(b) W^{-1}$. For the matrix units $\{ E_{ij} | 1 \leq i,j \leq n \}$ we have as above $W E_{ii} W^* = E_{(i-1) (i-1)}$ for $2 \leq i \leq n-1$ and $W E_{11} W^* = E_{nn}$. So  $$W \begin{pmatrix} \pi_1(b) & 0 & ... & 0  \\ 0 & \pi_2(b) & ... & 0 \\   . & . & . & . \\ 0 & 0 & ... & \pi_n(b) \end{pmatrix} W^* = \begin{pmatrix} \pi_2(b) & 0 & ... & 0  \\ 0 & \pi_3(b) & ... & 0 \\   . & . & . & . \\ 0 & 0 & ... & \pi_1(b) \end{pmatrix} ,$$
$$ = \begin{pmatrix} \pi_1(\Ad(u)(b)) & 0 & ... & 0  \\ 0 & \pi_2(\Ad(u)(b)) & ... & 0 \\   . & . & . & . \\ 0
& 0 & ... & \pi_n(\Ad(u)(b)) \end{pmatrix} = \pi(\Ad(u)(b)),$$ just what we wanted. 
\par
Now for $b \in B$ and $0 \leq k \leq n-1$ we have $$\tilde{\pi}((b u^k)^*) = \tilde{\pi}(u^{-k} b^*) = 
\tilde{\pi}(u^{-k} b^* u^k u^{-k}) = \pi(u^{-k} b^* u^k) W^{-k} = W^{-k} \pi(b^*) W^k W^{-k} , $$
$$ = W^{-k} \pi(b)^* = (\pi(b) W^k)^* = (\tilde{\pi}(b u^k))^*.$$ 
Also if $b, b' \in B$ and $0 \leq k, k' \leq n-1$, then $$\tilde{\pi}((b' u^{k'}).(b u^k)) = \tilde{\pi}(b'(u^{k'} b u^{-k'}) u^{k+k'}) = 
 \pi(b'(u^{k'} b u^{-k'})) W^{k+k'},$$ 
 $$= \pi(b') \pi(u^{k'} b u^{-k'}) W^{k+k'} = \pi(b') W^{k'} \pi(b) W^{-k'} W^{k+k'} =  
\tilde{\pi}(b' u^{k'}) \tilde{\pi}(b u^k).$$ This proves that that $\tilde{\pi}$ is a $*$-homomorphism. Continuity follows from continuity of $\pi$ and the Banach space representation $A = \underset{i=0}{\overset{n-1}{\bigoplus}} Bu^i$.
\par
Clearly $A_0 = \underset{i=0}{\overset{n-1}{\bigoplus}} B_0 u^i$ as a Banach space. 
It's also clear by the definition of $\tilde{\pi}$ 
that $A_0 \subset \ker(\tilde{\pi})$. Since $A_0$ has a Banach space codimension $n^2$ in $A$, and so does 
$\ker(\tilde{\pi})$, then we must have $A_0 = \ker(\tilde{\pi})$. 
\par
From the construction of the map $\tilde{\pi}$ we see that $\tilde{\pi}(e_{ii}) = E_{ii}$, since $\pi(e_{ii}) = E_{ii}$ and also 
$\tilde{\pi}(u^k) = W^k$. Since $\{e_{ii} | 1 \leq i \leq n \} \cup \{ W^k | 0 \leq k \leq n-1 \}$ generate $\mathbb{M}_n$, then we have 
$\tilde{\pi}(e_{ij}) = E_{ij}$, so the inclusion map $s: \mathbb{M}_n \rightarrow A$ given by $E_{ij} \mapsto e_{ij}$ is a right inverse
for  $\tilde{\pi}$.

\end{proof}

From this lemma follows that we can write $A = A_0 \oplus \mathbb{M}_n$ as a Banach space.

\begin{lemma}

If $\eta$ is a trace on $A_0$, then the linear functional on $A$ $\tilde{\eta}$, defined by 
$\tilde{\eta}(a_0 \oplus M) = \eta(a_0) + tr_n(M)$, where $a_0 \in A_0$ and $M \in \mathbb{M}_n$ 
is a trace and $\tilde{\eta}$ is the unique extension of $\eta$ to a trace on $A$ (of norm 1).

\end{lemma}

\begin{proof}

The functional $\eta$ can be extended in 
at most one way to a tracial state on $A$, because of the requirement $\tilde{\eta}(1_A) = 1$, the fact that $\mathbb{M}_n$ sits as a 
subalgebra in $A$, and the uniqueness on trace on $\mathbb{M}_n$. Since $\tilde{\eta}(1_A) = 1$, to show that $\tilde{\eta}$ is a trace we 
need to show that $\tilde{\eta}$ is positive and satisfies the trace property. For the trace property: If $x ,y \in A$ then we need to show
$\tilde{\eta}(xy) = \tilde{\eta}(yx)$. It is easy to see, that to prove this it's enough to prove that if $a_0 \in A_0$ and $M \in
\mathbb{M}_n$, then $\eta(a_0 M) = \eta(M a_0)$. Since $\eta$ is linear and $a_0$ is a linear combination of 4 positive elements we can
think, without loss of generality, that $a_0 \geq 0$. Then $a_0 = a_0^{1/2} a_0^{1/2}$ and $M a_0^{1/2}, a_0^{1/2} M \in A_0$, so since $\eta$ is a trace on $A_0$, we
have $\eta(M a_0) = \eta((M a_0^{1/2}) a_0^{1/2}) = \eta(a_0^{1/2}(M a_0^{1/2})) = \eta((a_0^{1/2} M) a_0^{1/2}) = \eta(a_0^{1/2}(a_0^{1/2}
M)) = \eta(a_0 M).$ This shows that $\tilde{\eta}$ satisfies the trace property. It remains to show positivity. Suppose $a_0 \oplus M \geq
0$. We must show $\eta(a_0 \oplus M) \geq 0$. Write $M = \underset{i=0}{\overset{n}{\sum}} \underset{j=0}{\overset{n}{\sum}} m_{ij} e_{ij}$ and
$a_0 = \underset{i=0}{\overset{n}{\sum}} \underset{j=0}{\overset{n}{\sum}} e_{ii} a_0 e_{jj}$ Since $\tilde{\eta}$ is a trace if $i \neq
j$, then $\tilde{\eta}(e_{ii} a_0 e_{jj}) = \tilde{\eta}(e_{jj} e_{ii} a_0) = 0$, so this shows that $\tilde{\eta}(a_0 \oplus M) = 
\underset{i=0}{\overset{n}{\sum}} (\frac{m_{ii}}{n} + \eta(e_{ii} a_0 e_{ii}))$. Clearly $a_0 \oplus M \geq 0$ implies $\forall 1 \leq i
\leq n, e_{ii} (a_0 \oplus M) e_{ii} \geq 0$. So to show positivity we only need to show $\forall 1 \leq i \leq n$ $\tilde{\eta}(e_{ii}(a_0
+ M)e_{ii}) \geq 0$, given $\forall 1 \leq i \leq n, m_{ii} e_{ii} + e_{ii} a_0 e_{ii} \geq 0$. Suppose that for some $i$, $m_{ii} < 0$. 
Then it follows that $e_{ii} a_0 e_{ii} \geq -m_{ii} e_{ii}$, so $e_{ii} a_0 e_{ii} \in e_{ii} A_0 e_{ii}$ is invertible, which implies 
$e_{ii} \in A_0$, that is not true. So this shows that $m_{ii} \geq 0$, and $m_{ii} e_{ii} \geq -e_{ii} a_0 e_{ii}$. If 
$\{ \epsilon_{\gamma} \}$ is an approximate unit for $A_0$, then positivity of $\eta$ implies 
$1 = \| \eta \| = \underset{\gamma}{\lim}\ 
\eta(\epsilon_{\gamma})$. Since $\eta$ is a trace we have 
$\underset{\gamma}{\lim}\ \eta(\epsilon_{\gamma} e_{ii}) = \frac{1}{n}$. Since 
$\forall \gamma,\ m_{ii}\epsilon_{\gamma}^{1/2} e_{ii} \epsilon_{\gamma}^{1/2} \geq - 
\epsilon_{\gamma}^{1/2} e_{ii} a_0 e_{ii} 
\epsilon_{\gamma}^{1/2}$, then $$tr_n(m_{ii} e_{ii}) = \frac{m_{ii}}{n} = 
\underset{\gamma}{\lim}\ \eta(m_{ii} e_{ii} \epsilon_{\gamma}) = 
\underset{\gamma}{\lim}\ \eta(m_{ii} \epsilon_{\gamma}^{1/2} e_{ii} \epsilon_{\gamma}^{1/2}) 
\geq \underset{\gamma}{\lim}\ 
\eta(\epsilon_{\gamma}^{1/2} e_{ii} a_0 e_{ii} \epsilon_{\gamma}^{1/2}),$$  
$$ = \underset{\gamma}{\lim}\ \eta(e_{ii} a_0 e_{ii} \epsilon_{\gamma}) 
= \eta(e_{ii} a_0 e_{ii}).$$ This finishes the proof of positivity and the proof of the lemma.

\end{proof}

\begin{remark}

We will show below that $\tau|_{A_0}$ is the unique trace on $A_0$. Since we have $A = A_0 \oplus \mathbb{M}_n$ as a Banach space, then 
clearly the free product trace $\tau$ on $A$ is given by $\tau(a_0 \oplus M) = \tau|_{A_0}(a_0) + tr_n(M)$, where $a_0 \oplus M \in A_0 
\oplus \mathbb{M}_n = A$. All tracial positive linear functionals of norm $\leq 1$ on $A_0$ are of the form $t\tau|_{A_0}$, where $0 \leq t \leq 1$. 
Then there will be no other traces on $A$ then the family $\lambda_t \overset{def}{=} t \tau|_{A_0} \oplus tr_n$. To show that these are 
traces indeed, we can use the above lemma (it is still true, no mater that the norm of $t \tau_{A_0}$ can be less than one), or we can 
represent them as a convex linear combination $\lambda_t = t \tau + (1-t)\mu$ of the free product trace $\tau$ and the trace $\mu$, 
defined by $\mu(a_0 \oplus M) = tr_n(M) = tr_n(\tilde{\pi}(a_0 \oplus M))$.

\end{remark}

\begin{lemma}

$\bar{B_0}$ (in $\mathfrak{H}_A$) $=$ $\bar{B}$ (in $\mathfrak{H}_A$).

\end{lemma}

\begin{proof}

 Let's take 
$D \overset{def}{=} ( \overset{1-\tilde{p}}{\mathbb{C}} \oplus \overset{\tilde{p}}{\mathbb{C}} ) * ( \overset{e_{11}}{\mathbb{C}} \oplus 
\overset{e_{22} + ... + e_{nn}}{\mathbb{C}}) \subset B$. Denote $D_0 \overset{def}{=} D\cap B_0$. From Theorem 2.3 follows that 
$D \cong \{ f: [0,b] \rightarrow \mathbb{M}_2 | f$ is continuous and $f(0)$ - diagonal$\}$ $\oplus \overset{\tilde{p}\wedge (1-e_{11})}
{\mathbb{C}}$, where $0 < b < 1$ and $\tau|_D$ is given by an atomless measure $\mu$ on $\{ f: [0,b] \rightarrow \mathbb{M}_2 | f$ is 
continuous and $f(0)$ - diagonal $\}$, $\tilde{p}$ is represented by $\begin{pmatrix} 1 & 0 \\ 0 & 0 \end{pmatrix} \oplus 1$, and $e_{11}$ 
is represented by $\begin{pmatrix} 1-t & \sqrt{t(1-t)} \\ \sqrt{t(1-t)} & t \end{pmatrix} \oplus 0$. A $*$-homomorphism, defined on the
generators of a $C^*$-algebra can be extended in at most one way to the whole $C^*$-algebra. This observation, together with $\pi_1(e_{11}) = 
\pi_1(\tilde{p}) = 1$ and $\pi_i(e_{22} = ... + e_{nn}) = \pi(\tilde{p}) = 1$ implies that $\pi_1|_D(f \oplus c) = f_{11}(0)$ and
$\pi_i|_D(f \oplus c) = c$ for $2 \leq i \leq n-1$. This means that $D_0 = \{ f: [0,b] \rightarrow \mathbb{M}_2 | f$ is continuous and 
$f_{11}(0) = f_{12}(0) = f_{21}(0) = 0 \} \oplus 0$. Now we see $\bar{D_0}$ (in $\mathfrak{H}_D$) $\cong$ $\mathbb{M}_2 
\otimes L^{\infty}([0,b], \mu) \oplus 0$, so then $e_{11} \in \bar{D_0}$ (in $\mathfrak{H}_D$). 
So we can find sequence $\{ \varepsilon_n \}$ of self-adjoined elements (functions) of $D_0$, supported on $e_{11}$, weakly converging to 
$e_{11}$ on $\mathfrak{H}_D$ and such that $\{ \varepsilon_n^2 \}$ also converges weakly to $e_{11}$ on $\mathfrak{H}_D$. 
Then take $a_1, a_2 \in A$. in $\mathfrak{H}_A$ we have $\langle \widehat{a_1}, (\varepsilon_n^2 - e_{11})\widehat{a_2} \rangle = 
\tau( (\varepsilon_n^2 - e_{11}) a_2 a_1^*) = \tau((\varepsilon_n - e_{11}) a_2 a_1^* (\varepsilon_n - e_{11})) \leq 4 \| a_2 a_1^* \| 
\tau(\varepsilon_n^2 - e_{11})$ (The last inequality is obtained by representing $a_2 a_1^*$ as a linear combination of 4 positive elements and using Cauchy-Bounjakovsky-Schwartz inequality). This shows that $e_{11} \in \bar{D_0}$ (in $\mathfrak{H}_A$) $\subset \bar{B_0}$ (in $\mathfrak{H}_A$). 
Analoguously $e_{ii} \in \bar{B_0}$ (in $\mathfrak{H}_A$), so this shows $\bar{B_0} = \bar{B}$ (in $\mathfrak{H}_A$). 

\end{proof}

It easily follows now that

\begin{cor}

$\bar{A_0}$ (in $\mathfrak{H}_A$) $=$ $\bar{A}$ (in $\mathfrak{H}_A$).

\end{cor}

The representation of $B_0$ on $\mathfrak{H}_A$ is faithful and nondegenerate, and we can use Proposition 3.6, together with Theorem 3.5 and
the fact that $\Ad(u^k)$ are outer for $\bar{B} = \bar{B_0}$ to get:

\begin{lemma}

$A_0 = B_0 \rtimes G$ is simple.

\end{lemma} 

For the uniqueness of trace we need to modify a little the proof Theorem 3.7 (which is Theorem 1 in \cite{B93}, stated for "nontwisted"
crossed products).

\begin{lemma}

$A_0 = B_0 \rtimes G$ has a unique trace, $\tau|_{A_0}$.

\end{lemma}

\begin{proof}

Above we already proved that $\{ \Ad(u^k) | 1 \leq k \leq n-1 \}$ are 
$\tau|_{B_0}$-outer for $B_0$. 
\par

Suppose that $\eta$ is a trace on $A_0$. We will show that $\tau|_{A_0} = \eta$. We consider the GNS 
representation of $B$, associated 
to $\tau|_B$. By repeating the proof of Lemma 3.13 we see that $\bar{B_0}$ (in $\mathfrak{H}_B$) $=$ 
$\bar{B}$ (in $\mathfrak{H}_B$). 
The simplicity of $B_0$ 
allows us to identify $B_0$ with $\pi_{\tau|_B}(B_0)$. We will also identify $B_0$ with it's canonical copy in $A_0$. $A_0$ is 
generated by $\{ b_0 \in B_0 \} \cup \{ u^k | 0 \leq k \leq n-1 \}$ and $\{ \Ad(u^k) | 0 \leq k \leq n-1 \}$ 
extend to $\bar{B_0}$ 
(in $\mathfrak{H}_A$), so also to $\bar{B_0}$ (in $\mathfrak{H}_B$) ( $\cong \bar{B}$ (in $\mathfrak{H}_A$)). 
Now we can form the von Neumann algebra crossed product $\tilde{A} \overset{def}{=} \bar{B_0} \rtimes 
\{ \Ad(u^k) | 0 \leq k \leq n-1 \} \cong \bar{B} \rtimes \{ \Ad(u^k) | 0 \leq k \leq n-1 \}$, where the weak closures are in 
$\mathfrak{H}_B$. Clearly $\tilde{A} \cong \bar{A}$ (in $\mathfrak{H}_A$). Denote by $\widetilde{\tau_{B_0}}$ the extension of $\tau|_{B_0}$ to 
$\bar{B_0}$ (in $\mathfrak{H}_A$), given by $\widetilde{\tau_{B_0}}(x) = \langle x(\widehat{1_A}), \widehat{1_A} \rangle_{\mathfrak{H}_A}$. By Proposition 3.19 of 
Chapter V in \cite{T79}, $\widetilde{\tau_{B_0}}$ is a faithful normal trace on $\bar{B_0}$ (in $\mathfrak{H}_A$). Now from the fact that 
$\bar{B_0}$ (in $\mathfrak{H}_A$) is a factor and using Lemma 1 from \cite{L81} we get that $\widetilde{\tau_{B_0}}$ is unique on 
$\bar{B_0}$ (in $\mathfrak{H}_A$). By the same argument we have that the extension $\widetilde{\tau_{A_0}}$ of $\tau|_{A_0}$ to $\bar{A_0}$ 
(in $\mathfrak{H}_{A}$) $\cong$ $\bar{A}$ (in $\mathfrak{H}_A$) is unique, since $\bar{A_0}$ (in $\mathfrak{H}_{A}$) $\cong$ $\bar{A}$ 
(in $\mathfrak{H}_A$) is a factor. 
\par
We take the unique extension of $\eta$ to $A$. We will call it again $\eta$ for convenience. \\ 
We denote by $\mathfrak{H}'_{C}$ the GNS Hilbert space for $C$, corresponding to $\eta|_C$ (for $C$ $=$ $A$, $B$, $B_0$, $A_0$). 
Since $\eta|_{B_0} = \tau|_{B_0}$ it follows that $\bar{B_0}$ (in $\mathfrak{H}'_{B_0}$) $\cong$ 
$\bar{B}$ (in $\mathfrak{H}'_B$) and of course $\mathfrak{H}'_{B_0} = \mathfrak{H}'_B$. Then similarly as in Lemma 3.12 we get that 
$\bar{A_0}$ (in $\mathfrak{H}'_{A_0}$) $\cong$ $\bar{A}$ (in $\mathfrak{H}'_{A}$), so 
$\mathfrak{H}'_{A_0} = \mathfrak{H}'_{A}$ (this can be done, since $\tau|_{B_0} = \eta|_{B_0}$). 
Now again by 
Proposition 3.19 of Chapter V in \cite{T79} we have that 
$\tilde{\eta}(x) \overset{def}{=} \langle x(\widehat{1_A}), \widehat{1_A} \rangle_{\mathfrak{H}'_A}$ 
($\widehat{1_A}$ is abuse of notation - in this case it's the element, corresponding to $1_A$ in
$\mathfrak{H}'_A$) defines a 
faithful normal trace on $\overline{\pi'_{A}(A)}$ (in $\mathfrak{H}'_A$). 
In particular $\tilde{\eta}|_{\overline{\pi'_A(B)}}$ is a faithful normal trace on 
$\overline{\pi'_A(B)}$ (in $\mathfrak{H}'_{A}$). 
By uniqueness of $\tau|_{B_0}$ we have $\tau|_{B_0} = \eta|_{B_0}$, so for $b_0 \in B_0$ we have $\tilde{\tau}
(b_0) = \tau(b_0) = \eta(b_0) = \langle \pi'_{A}(b_0)(\widehat{1_A}), \widehat{1_A} \rangle_{\mathfrak{H}'_A} = \tilde{\eta}(\pi'_{A}(b_0))$. 
\par
Since $B_0$ is simple, it follows that $\pi'_{A}|_{B_0}$ is a $*$-isomorphism from $B_0$ onto $\pi'_{A}(B_0)$ and from Exercise 7.6.7 in 
\cite{KR86} it follows that $\pi'_{A}|_{B_0}$ extends to a $*$-isomorphism from $\bar{B_0}$ (in $\mathfrak{H}_A$) $\cong$ $\bar{B}$ 
(in $\mathfrak{H}_A$) onto $\overline{\pi'_{A}(B_0)}$ (in $\mathfrak{H}'_A$) $\cong$ $\overline{\pi'_{A}(B)}$ (in $\mathfrak{H}'_A$). 
We will denote this $*$-isomorphism by $\theta$. 
We set $w \overset{def}{=} \pi'_A(u)$, $\beta \overset{def}{=} \theta \Ad(u) \theta^{-1} 
\in Aut(\overline{\pi'_A(B)}$ (in $\mathfrak{H}'_A$)). For $b_0 \in B_0$ we have 
$w \pi'_A(b_0) w^* = \pi'_A(u b_0 u^*) = \pi'_A((\Ad(u))(b_0)) = 
\beta (\pi'_A(b_0))$. 
So by weak continuity follows $\beta = \Ad(w)$ on $\overline{\pi'_A(B)}$ (in $\mathfrak{H}'_A$). 
Since $\bar{B}$ 
(in $\mathfrak{H}_A$) is a factor and $\{ \Ad(u^k) | 1 \leq k \leq n-1 \}$ are all outer, Kallman's Theorem 
(Corrolary 1.2 in \cite{Ka69}) gives us that $\{ \Ad(u^k) | 1 \leq k \leq n-1 \}$ act freely on $\bar{B}$ 
(in $\mathfrak{H}_A$). Namely if $\bar{b} \in \bar{B}$ (in $\mathfrak{H}_A$), and if $\forall \bar{b}'
\in \bar{B}$ (in $\mathfrak{H}_A$), $\bar{b} \bar{b}' = \Ad(u^k)(\bar{b}') \bar{b}$, then $\bar{b} = 0$. 
Then by the above settings it is
clear that $\{ \Ad(w^k) | 1 \leq k  \leq n-1 \}$ also act freely on $\overline{\pi'_A(B)}$ 
(in $\mathfrak{H}'_A$). 
\par
Since $\tilde{\eta}$ is a faithful normal trace on $\overline{\pi'_A(A)}$ (in $\mathfrak{H}'_A$), then by Proposition 2.36 of Chapter V in
\cite{T79} there exists a faithful conditional expectation $P: \overline{\pi'_A(A)} \rightarrow \overline{\pi'_A(B)}$ (both weak 
closures are in $\mathfrak{H}'_A$). 
$\forall x \in \overline{\pi'_A(B)}$ (in $\mathfrak{H}'_A$), and $\forall 1 \leq k \leq n-1$, 
$\Ad(w^k)(x) w^k = w^k x$. Applying $P$ we get $\Ad(w^k)(x)(P(w^k)) = P(w^k) x$, so by the free action of 
$\Ad(w^k)$ we get that $P(w^k) = 0$, $\forall 1 
\leq k \leq n-1$. 
It's clear that $\{ \overline{\pi'_A(B)} \} \cup \{ w^k | 1 \leq k  \leq n-1 \}$ generates 
$\overline{\pi'_A(A)}$ (in
$\mathfrak{H}'_A$) as a von Neumann algebra. 
Now we use Proposition 22.2 from \cite{S81}. It gives us a $*$-isomorphism $\Phi : 
\overline{\pi'_A(A)}$ (in $\mathfrak{H}'_A$) $\rightarrow \bar{B} \rtimes 
\{ \Ad(u^k) | 1 \leq k \leq n-1 \} \cong \bar{A}$ (last two weak 
closures are in $\mathfrak{H}_A$) with $\Phi(\theta(x)) = x,$ $x\in \bar{B}$ (in $\mathfrak{H}_A$), 
$\Phi(w) = u$. So since $\bar{A}$ 
(in $\mathfrak{H}_A$) is a finite factor, so is $\overline{\pi'_A(A)}$ (in $\mathfrak{H}'_A$), and so it's trace $\tilde{\eta}$ is unique. 
Hence, $\tilde{\eta} = \tilde{\tau} \circ \Phi$, and so $\forall b \in B$, and $\forall 1 \leq k \leq n-1$ we have $\eta(b u^k) = 
\tilde{\eta}(\pi'_A(b) \pi'_A(u^k)) = \tilde{\tau}(\Phi(\pi'_A(b)) \Phi(\pi'_A(u^k))) = \tilde{\tau}(\Phi(\theta(b)) \Phi(w^k)) = 
\tilde{\tau}(b u^k) = \tau(b u^k)$. By continuity and linearity of both traces we get $\eta = \tau$, just what we want.

\end{proof}

We conclude this section by proving the following 

\begin{prop}

Let 

\begin{center}

$(A,\tau) \overset{def}{=} ( \underset{\alpha_1}{\overset{p_1}{\mathbb{C}}} \oplus ... \oplus \underset{\alpha_m}{\overset{p_m}
{\mathbb{C}}})*(\mathbb{M}_n, tr_n)$, 

\end{center}

where $\alpha_1 \leq \alpha_2 \leq ... \leq \alpha_m$. Then: 
\par
(I) If $\alpha_m < 1-\frac{1}{n^2}$, then $A$ is unital, simple with a unique trace $\tau$. 
\par
(II) If $\alpha_m = 1-\frac{1}{n^2}$, then we have a short exact sequence $0 \rightarrow A_0 \rightarrow A \rightarrow \mathbb{M}_n
\rightarrow 0$, where $A$ has no central projections, and $A_0$ is nonunital, simple with a unique trace $\tau|_{A_0}$. 
\par
(III) If $\alpha_m > 1-\frac{1}{n^2}$, then $A = \underset{n^2 - n^2 \alpha_m}{\overset{f}{A_0}} \oplus \underset{n^2 \alpha_m - n^2 +
1}{\overset{1-f}{\mathbb{M}_n}}$, where $1-f \leq p_m$, and where $A_0$ is unital, simple and has a unique trace 
$(n^2 - n^2 \alpha_m)^{-1} \tau|_{A_0}$. 
\par
Let $f$ means the identity projection for cases (I) and (II). Then in all cases for each of the projections $f p_1, ..., f p_m$ we have a unital, 
diffuse abelian $C^*$-subalgebra of $A$, supported on it.

\par
In all the cases $p_m$ is a full projection in $A$.

\end{prop}

\begin{proof}

We have to prove the second part of the proposition, since the first part follows from Lemma 3.10, Lemma 3.11, Lemma 3.12, Lemma 3.17 and
Lemma 3.18. From the discussion above we see that in all cases we have $fA = fB \rtimes \{ \Ad(f u^k f) | 0 \leq k \leq n-1 \}$, where $B$
and $\{ \Ad(f u^k ) | 0 \leq k \leq n-1 \}$ are as above. So the existence of the unital, diffuse abelian $C^*$-sublagebras follows from 
Theorem 2.4, applied to $B$. 
\par
In the case (I) $p_m$ is clearly full, since $A$ is simple. In the case (III) it's easy to see that $p_m \wedge f \neq 0$ and $p_m \geq (1-f)$, so since $A_0$ and $\mathbb{M}_n$ are simple in this case, then $p_m$ is full in $A$. In case (II) it follows from Theorem 2.4 that $p_m$ is full in $B$, and consequently in $A$.

\end{proof}

\section{ The General Case}

In this section we prove the general case of Theorem 2.6, using the result from the previous section
(Proposition 3.19). The prove of the general case involves techniques from \cite{D99}. So we will need two technical results from there.
\par 
The first one is Proposition 2.8 in \cite{D99} (see also \cite{D93}):

\begin{prop}

Let $A = A_1 \oplus A_2$ be a direct sum of unital $C^*$-algebras and
let $p = 1 \oplus 0 \in A$. Suppose $\phi_A$ is a state on $A$ with $0
< \alpha \overset{def}{=} \phi_A(p) < 1$. Let $B$ be a unital
$C^*$-algebra with a state $\phi_B$ and let $(\mathfrak{A}, \phi) =
(A, \phi_A) * (B, \phi_B)$. Let $\mathfrak{A}_1$ be the
$C^*$-subalgebra of $\mathfrak{A}$ generated by $(0 \oplus A_2) +
\mathbb{C} p \subseteq A$, toghether with $B$. In other words 

\begin{equation*}
(\mathfrak{A}_1, \phi|_{\mathfrak{A}_1}) =
(\underset{\alpha}{\overset{p}{\mathbb{C}}} \oplus 
\underset{1-\alpha}{\overset{1-p}{A_2}}) * (B, \phi_B).
\end{equation*}

Then $p \mathfrak{A} p$ is generated by $p \mathfrak{A}_1 p$ and $A_1
\oplus 0 \subset A$, which are free in $(p \mathfrak{A} p,
\frac{1}{\alpha} \phi|_{p \mathfrak{A} p})$. In other words 

\begin{equation*}
(p \mathfrak{A} p, \frac{1}{\alpha} \phi|_{p \mathfrak{A} p}) \cong 
(p \mathfrak{A}_1 p,\frac{1}{\alpha} \phi|_{p \mathfrak{A}_1 p}) *
(A_1, \frac{1}{\alpha} \phi_A|_{A_1}).
\end{equation*}

\end{prop}

\begin{remark}

This proposition was proved for the case of von Neumann algebras in \cite{D93}. It is true also in the case of $C^*$-algebras.

\end{remark}

The second result is Proposition 2.5 (ii) of \cite{D99}, which is easy
and we give its proof also:

\begin{prop}

Let $A$ be a $C^*$-algebra. Take $h \in A, h \geq 0$, and let $B$ be the hereditary subalgebra 
$\overline{hAh}$ of $A$ ( $\overline{*}$ means norm closure). Suppose that $B$ is full in $A$. Then if 
$B$ has a unique trace, then $A$ has at most one tracial state.

\end{prop}

\begin{proof}

It's easy to see that $\Span \{ xhahy | a,x,y \in A \}$ is norm dense
in $A$. If $\tau$ is a tracial state on $A$ then $\tau(xhahy) =
\tau(h^{1/2} ahyx h^{1/2})$. Since $h^{1/2} ahyx h^{1/2} \in B$,
$\tau$ is uniquely determined by $\tau_B$.

\end{proof}

It is clear that Proposition 3.19 agrees with Theorem 2.6, so it is 
a special case.
\par
As a next step we look at a $C^*$-algebra of the form 

\begin{equation*}
(M, \tau) = (\underset{\alpha_0'}{\overset{p_0'}{A_0}} \oplus \underset{\alpha_1'}{\overset{p_1'}{\mathbb{M}_{m_1}}} \oplus ... \oplus \underset{\alpha_k'}{\overset{p_k'}{\mathbb{M}_{m_k}}} \oplus \underset{\alpha_1}{\overset{p_1}{\mathbb{C}}} \oplus ... \oplus \underset{\alpha_l}{\overset{p_l}{\mathbb{C}}}) * (\mathbb{M}_n, tr_n),
\end{equation*}
where $A_0$ comes with a specified trace and has a unital, diffuse
abelian $C^*$-subalgebra with unit $p'_0$. Also we suppose that
$\alpha'_0 \geq 0$, $0 < \alpha_1' \leq ... \leq \alpha_k'$, $0 <
\alpha_1 \leq ... \leq \alpha_l$, $m_1, ..., m_k \geq 2$, and either $\alpha'_0 > 0$ or $k \geq 1$, or both. Let's denote $p_0 
\overset{def}{=} p_0' + p_1' + ... + p_k'$, $B_0 \overset{def}{=} 
\underset{\alpha_1'}{\overset{p_1'}{\mathbb{M}_{m_1}}} \oplus ... 
\oplus \underset{\alpha_k'}{\overset{p_k'}{\mathbb{M}_{m_k}}}$, and $\alpha_0 \overset{def}{=} \alpha'_0 + \alpha'_1 + ... + \alpha'_k = \tau(p_0)$. \\ 
Let's have a look at the $C^*$-subalgebras $N$  and $N'$ of $M$ given by 

\begin{equation*}
(N, \tau|_N) = (\underset{\alpha_0}{\overset{p_0}{\mathbb{C}}} \oplus \underset{\alpha_1}{\overset{p_1}{\mathbb{C}}} \oplus ... \oplus \underset{\alpha_l}{\overset{p_l}{\mathbb{C}}}) * (\mathbb{M}_n, tr_n)
\end{equation*}
and

\begin{equation*}
(N', \tau|_{N'}) = (\underset{\alpha_0'}{\overset{p_0'}{\mathbb{C}}} \oplus \underset{\alpha_1'}{\overset{p_1'}{\mathbb{C}}} \oplus ... \oplus \underset{\alpha_k'}{\overset{p_k'}{\mathbb{C}}} \oplus \underset{\alpha_1}{\overset{p_1}{\mathbb{C}}} \oplus ... \oplus \underset{\alpha_l}{\overset{p_l}{\mathbb{C}}}) * (\mathbb{M}_n, tr_n).
\end{equation*}

We studied the $C^*$-algebras, having the form of $N$ and $N'$ in the previous section. A brief description is as follows:
\par
If $\alpha_0, \alpha_l < 1-\frac{1}{n^2}$, then $N$ is
simple with a unique trace and $N'$ is also simple with a unique trace. For each of the projections $p_0', p_1', ..., p_k', p_1, ..., p_l$ we have a unital, diffuse abelian $C^*$-subalgebra of $N'$, supported on it. 
\par 
If $\alpha_0$, or $\alpha_l$ $= 1-\frac{1}{n^2}$, then
$N$ has no central projections, and we have a short exact sequence 
$0 \rightarrow N_0 \rightarrow N \rightarrow \mathbb{M}_n \rightarrow 
0$, with $N_0$ being simple with a unique trace. Moreover $p_0$ or $p_l$ respectivelly is full in $N$. For each of the
projections $p_0', p_1', ..., p_k', p_1, ..., p_l$ we have a unital,
diffuse abelian $C^*$-subalgebra of $N'$, supported on it. 
\par
If $\alpha_0$ or $\alpha_l$ $> 1-\frac{1}{n^2}$, then
$N = \overset{q}{N_0} \oplus \mathbb{M}_n$, with $N_0$ being simple
and having a unique trace. 
\par
We consider 2 cases: 
\par
(I) case:  $\alpha_l \geq \alpha_0$. 
\par
(1) $\alpha_l < 1-\frac{1}{n^2}$. 
\par
In this case $N$ and $N'$ are simple and has unique traces, and $p_0$ is full in $N$ and consequently $1_M = 1_N$ is contained in $\langle p_0 \rangle_N$ - the ideal of $N$, generated by $p_0$. Since $\langle p_0 \rangle_N \subset \langle p_0 \rangle_M$ it follows that $p_0$ is full also in $M$. From Proposition
4.1 we get $p_0 M p_0 \cong (A_0 \oplus B_0) * p_0 N p_0$. Then from
Theorem 3.9 follows that $p_0 M p_0$ is simple and has a unique trace.
Since $p_0$ is a full projection, Proposition 4.3 tells us that 
$M$ is simple and $\tau$ is its unique trace. For each of the
projections $p_0', p_1', ..., p_k', p_1, ..., p_l$ we have a unital,
diffuse abelian $C^*$-subalgebra of $M$, supported on it, and comming 
from $N'$.
\par
(2) $\alpha_l = 1-\frac{1}{n^2}$.
\par
In this case it is also true that for each of the projections 
$p_0', p_1', ..., p_k', p_1, ..., p_l$ we have a unital, diffuse 
abelian $C^*$-subalgebra of $M$, supported on it, and comming from 
$N'$. It is easy to see that $M$ is the linear span of $p_0 M p_0$, 
$p_0 M (1- p_0) N (1- p_0)$, $(1 -p_0) N p_0 M p_0$, $(1- p_0) N p_0 M
p_0 N (1- p_0)$ and $(1- p_0) N (1- p_0)$. We know that we have a
$*$-homomorphism $\pi : N \rightarrow M_n$, such that $\pi(p_l) = 1$.
Then it is clear that $\pi(p_0) = 0$, so we can extend $\pi$ to a
linear map $\tilde{\pi}$ on $M$, defining it to equal $0$ on $p_0 M 
p_0$, $p_0 M (1- p_0) N (1- p_0)$, $(1 -p_0) N p_0 M p_0$ and $(1- p_0) N 
p_0 M p_0 N (1- p_0)$. It is also clear then that $\tilde{\pi}$ will
actually be a $*$-homomorphism. Since $\ker(\pi)$ is simple in $N$ and
$p_0 \in \ker(\pi)$, then $p_0$ is full in $\ker(\pi) \subset N$, so by
the above representation of $M$ as a linear span we see that $p_0$ is
full in $\ker(\tilde{\pi})$ also. From Proposition 4.1 follows that
$p_0 M p_0 \cong (A_0 \oplus B_0) * (p_0 N p_0)$. Since $p_0 N p_0$
has a unital, diffuse abelian $C^*$-subalgebra with unit $p_0$, it
follows from Theorem 3.9 that $p_0 M p_0$ is simple and has a unique
trace (to make this conclusion we could use Theorem 1.5 instead). Now
since $p_0 M p_0$ is full and hereditary in $\ker(\tilde{\pi})$, from 
Proposition 4.3 follows that $\ker(\tilde{\pi})$ is simple and has a
unique trace.
\par
(3) $\alpha_l > 1-\frac{1}{n^2}$.
\par
In this case $N = \underset{n^2 - n^2 \alpha_l}{\overset{q}{N_0}} \oplus \underset{n^2 \alpha_l - n^2 + 1}{\overset{1-q}{\mathbb{M}_n}}$ and also $N' = \underset{n^2 - n^2 \alpha_l}{\overset{q}{N'_0}} \oplus \underset{n^2 \alpha_l - n^2 + 1}{\overset{1-q}{\mathbb{M}_n}}$ with $N_0$ and $N'_0$ being simple with unique traces. For each of the projections $q p_0', q p_1'
, ..., q p_k', q p_1, ..., q p_l$ we have a unital, diffuse abelian 
$C^*$-subalgebra of $M$, supported on it, and coming from $N'_0$.
\par 
Since $p_0 \leq q$ we can write $M$ as a linear span of $p_0 M p_0$, 
$p_0 M p_0 N_0 (1- p_0)$, $(1- p_0) N_0 p_0 M p_0$, $(1- p_0) N_0 p_0
M p_0 N_0 (1- p_0)$, $(1- p_0) N_0 (1- p_0)$ and $\mathbb{M}_n$. So we
can write $M = \underset{n^2 - n^2 \alpha_l}{\overset{q}{M_0}} \oplus \underset{n^2 \alpha_l - n^2 + 1}{\overset{1-q}{\mathbb{M}_n}}$,
where $M_0 \overset{def}{=} q M q \supset N_0$. 
We know that $p_0$ is full in
$N_0$, so  as before we can write $1_{M_0} = 1_{N_0} \in \langle p_0 \rangle_{N_0} \subset \langle p_0 \rangle_{M_0}$, so $\langle p_0 \rangle_{M_0} = M_0$. Because of Proposition 4.1, we can write $p_0 M_0 p_0 \cong (A_0
\oplus B_0) * (p_0 N_0 p_0)$. Since $p_0 N_0 p_0$ has a unital,
diffuse abelian $C^*$-subalgebra with unit $p_0$, then from Theorem
3.9 (or from Theorem 1.5) it follows that $p_0 M_0 p_0$ is simple with a
unique trace. Since $p_0 M_0 p_0$ is full and hereditary in $M_0$,
Proposition 4.3 yields that $M_0$ is simple with a unique trace.
\par
(II) $\alpha_0$ $>$ $\alpha_l$.
\par
(1) $\alpha_0 \leq 1- \frac{1}{n^2}$.
\par
In this case $p_0$ is full in $N$ and also in $N'$, so $1_M = 1_N \in \langle p_0 \rangle_N$, 
which means $p_0$ is full in $M$ also. $p_0 M p_0$ is a full hereditary $C^*$-subalgebra of $M$ and $p_0 M p_0 
\cong (A_0 \oplus B_0) * p_0 N p_0$ by Proposition 4.1. Since $p_0 N p_0$ has a diffuse abelian $C^*$-subalgebra, Theorem 3.9 
(or Theorem 1.5) shows that $p_0 M p_0$ is simple with a unique trace and then by Proposition 4.3 follows that the same is true for $M$.
For each of the projections $p_0', p_1', ..., p_k', p_1, ..., p_l$ we have a unital, diffuse abelian $C^*$-subalgebra of $M$, supported 
on it, comming from $N'$. 
\par
(2) $\alpha_0$ $> 1-\frac{1}{n^2}$. \\ 
We have 3 cases:
\par
(2$'$) $\alpha'_0 > 1-\frac{1}{n^2}$.
\par
In this case $N \cong \overset{q}{N_0} \oplus \mathbb{M}_n$ and $N' \cong \overset{q'}{N'_0} \oplus \mathbb{M}_n$, where $q \leq q'$, with $N_0$ and $N'_0$ being simple and having unique traces. It is easy to see that $p'_1, ..., p'_k, p_1, ..., p_l \leq q'$, so for each of the
projections $p_1', ..., p_k', p_1, ..., p_l$ we have a unital, diffuse abelian $C^*$-subalgebra of $N'$, supported on it. So those
$C^*$-subalgebras live in $M$ also. We have a unital, diffuse abelian $C^*$-subalgebra of $A_0$, supported on $1_{A_0}$, which yields a 
unital, diffuse abelian $C^*$-subalgebra on $M$, supported on $p'_0$. It is  clear that $p_0$ is full in $N$, so as before, $1_M =
1_N \in \langle p_0 \rangle_N$, so $p_0$ is full in $M$ also, so $p_0 M p_0$ is a full hereditary $C^*$-subalgebra of $M$. From
Proposition 4.1 we have $p_0 M p_0 \cong (A_0 \oplus B_0) * ( p_0 N_0 p_0 \oplus \mathbb{M}_n)$. It is easy to see that 
$\mathbb{M}_n$, for $n \geq 2$ contains two $tr_n$-orthogonal zero-trace unitaries. Since also $p_0 N_0 p_0$ has a 
unital, diffuse abelian $C^*$-subalgebra, supported on $1_{N_0}$, it is easy to see (using Proposition 2.2) that it also contains two $\tau|{N_0}$-orthogonal, zero-trace unitaries. Then the conditions of 
Theorem 1.5 are satisfied. This means that $p_0 M p_0$ is simple with a unique trace and Proposition 4.3 implies that $M$ is simple with 
a unique trace also.
\par
(2$''$) $\alpha'_k > 1-\frac{1}{n^2}$.
\par
Let's denote $$N'' = (\underset{\alpha_0'}{\overset{p_0'}{A_0}} \oplus \underset{\alpha_1'}{\overset{p_1'}{\mathbb{M}_{m_1}}} \oplus ... 
\oplus \underset{\alpha_{k-1}'}{\overset{p_{k-1}'}{\mathbb{M}_{m_{k-1}}}} \oplus \underset{\alpha_{k}'}{\overset{p_k'}{\mathbb{C}}} \oplus 
\underset{\alpha_1}{\overset{p_1}{\mathbb{C}}} \oplus ... \oplus \underset{\alpha_l}{\overset{p_l}{\mathbb{C}}}) * (\mathbb{M}_n, tr_n).$$ 
Then $N''$ satisfies the conditions of case (I,3) and so $N'' \cong \overset{q}{N''_0} \oplus \mathbb{M}_n$. Clearly $p_0', p_1', ...,
p_{k-1}', p_1, ..., p_l \leq q$, so for each of the projections $p_0', p_1', ..., p_{k-1}', p_1, ..., p_l$ we have a unital, diffuse 
abelian $C^*$-subalgebra of $N''_0$, supported on it. Those $C^*$-algebras live in $M$ also. From case (I,3) we have that $p'_k$ is full in 
$N''$ and as before $1_M = 1_{N''} \in \langle p'_k \rangle_{N''}$ implies that $p'_k$ is full in $M$ also. From Proposition 4.1 follows 
that $p'_k M p'_k \cong (p'_k N''_0 p'_k \oplus \mathbb{M}_n) * \mathbb{M}_{m_k}$. Since $N''_0$ has a unital, diffuse abelian 
$C^*$-subalgebra, supported on $q p'_k$, then an argument, similar to the one we made in case (II, 2$"$), allows to apply Theorem 1.5 to get that 
$p'_k M p'_k$ is simple with a unique trac. By Proposition 4.3 follows that the same is true for $M$. The unital, diffuse abelian 
$C^*$-subalgebra of $M$, supported on $p'_k$, we can get by applying the note after Theorem 1.5 to $p'_k M p'_k \cong 
(p'_k N''_0 p'_k \oplus \mathbb{M}_n) * \mathbb{M}_{m_k}$. 
\par
(2$'''$) $\alpha'_0$ and $\alpha'_k$ $\leq 1-\frac{1}{n^2}$. 
\par
In this case $N \cong \overset{q}{N_0} \oplus \mathbb{M}_n$, with $N_0$  being simple and having a unique trace. Moreover $N'$ has no central projections and for each of the projections $p'_0, p_1', ..., p_k', p_1, ..., p_l$ we have a unital, diffuse abelian $C^*$-subalgebra of $N'$, supported on it. So those
$C^*$-subalgebras live in $M$ also. It is  clear that $p_0$ is full in $N$, so as before $1_M = 1_N \in \langle p_0 \rangle_N$, so $p_0$ is full in $M$ also, so $p_0 M p_0$ is a full hereditary $C^*$-subalgebra of $M$. From
Proposition 4.1 we have $p_0 M p_0 \cong (A_0 \oplus B_0) * ( p_0 N_0 p_0 \oplus \mathbb{M}_n)$. Since $A_0$ and $p_0 N_0 p_0$ both have 
unital, diffuse abelian $C^*$-subalgebras, supported on their units, it is easy to see (using Proposition 2.2), that the conditions of 
Theorem 1.5 are satisfied. This means that $p_0 M p_0$ is simple with a unique trace and Proposition 4.3 yields that $M$ is simple with 
a unique trace also.
\par
We summarize the discussion above in the following 

\begin{prop}

Let 

\begin{equation*}
(M,\tau) \overset{def}{=} (\underset{\alpha_0'}{\overset{p_0'}{A_0}} \oplus \underset{\alpha_1'}{\overset{p_1'}{\mathbb{M}_{m_1}}} \oplus
... \oplus \underset{\alpha_k'}{\overset{p_k'}{\mathbb{M}_{m_k}}} \oplus \underset{\alpha_1}{\overset{p_1}{\mathbb{C}}} \oplus ... \oplus
\underset{\alpha_l}{\overset{p_l}{\mathbb{C}}}) * (\mathbb{M}_n, tr_n),
\end{equation*}
where $n \geq 2$, $\alpha'_0 \geq 0$, $\alpha'_1 \leq \alpha'_2 \leq ... \leq \alpha'_k$, $\alpha_1 \leq ... \leq \alpha_l$, $m_1, ..., m_k \geq 2$,
and $\overset{p'_0}{A_0} \oplus 0$ has a unital, diffuse abelian $C^*$-subalgebra, having $p'_0$ as a unit. Then: 
\par
(I) If $\alpha_l < 1-\frac{1}{n^2}$, then $M$ is unital, simple with a unique trace $\tau$. 
\par
(II) If $\alpha_l = 1-\frac{1}{n^2}$, then we have a short exact sequence $0 \rightarrow M_0 \rightarrow M \rightarrow \mathbb{M}_n
\rightarrow 0$, where $M$ has no central projections and $M_0$ is nonunital, simple with a unique trace $\tau|_{M_0}$. 
\par
(III) If $\alpha_l > 1-\frac{1}{n^2}$, then $M = \underset{n^2 - n^2 \alpha_l}{\overset{f}{M_0}} \oplus \underset{n^2 \alpha_l - n^2 +
1}{\overset{1-f}{\mathbb{M}_n}}$, where $1-f \leq p_l$, and where $M_0$ is unital, simple and has a unique trace 
$(n^2 - n^2 \alpha_l)^{-1} \tau|_{M_0}$. 
\par
Let $f$ means the identity projection for cases (I) and (II). Then in all cases for each of the projections $f p_0', f p_1', ..., f p_k', f p_1, ..., 
f p_l$ we have a unital, diffuse abelian $C^*$-subalgebra of $M$, supported on it.
\par
In all the cases $p_l$ is a full projection in $M$.

\end{prop}

To prove Theorem 2.6 we will use Proposition 4.4. First let's check that Proposition 4.4 agrees with the conclusion of
Theorem 2.6. We can write $$(M,\tau) \overset{def}{=} (\underset{\alpha_0'}{\overset{p_0'}{A_0}} \oplus \underset{\alpha_1'}{\overset{p_1'}
{\mathbb{M}_{m_1}}} \oplus ... \oplus \underset{\alpha_k'}{\overset{p_k'}{\mathbb{M}_{m_k}}} \oplus \underset{\alpha_1}{\overset{p_1}
{\mathbb{C}}} \oplus ... \oplus \underset{\alpha_l}{\overset{p_l}{\mathbb{C}}}) * \underset{\beta_1}{\overset{q_1}{\mathbb{M}_n}},$$ where
$q_1 = 1_M$ and $\beta_1 = 1$. It is easy to see that $L_0 = \{ (l,1) | \frac{\alpha_l}{1^2} + \frac{1}{n^2} = 1 \} = \{ (l,1) | \alpha_l =
1-\frac{1}{n^2} \}$, which is not empty if and only if $\alpha_l = 1-\frac{1}{n^2}$. Also $L_+ = \{ (l,1) | \frac{\alpha_l}{1^2} +
\frac{1}{n^2} > 1 \} = \{ (l,1) | \alpha_l > 1-\frac{1}{n^2} \}$, and here  $L_+$ is not empty if and only if 
$\alpha_l > 1-\frac{1}{n^2}$. If both $L_+$ and $L_0$ are empty, then $M$ is simple with a unique trace. If $L_0$ is not empty, then
clearly $L_+$ is empty, so we have no central projections and a short exact sequence $0 \rightarrow M_0 \rightarrow M \rightarrow
\mathbb{M}_n \rightarrow 0$, with $M_0$ being simple with a unique trace. In this case all nontrivial projections are full in $M$. If
$L_+$ is not empty, then clearly $L_0$ is empty and so $M = \underset{n^2 -n^2 \alpha_l}{\overset{q}{M_0}} \oplus 
\underset{n^2(\frac{\alpha_l}{1^2} + \frac{1}{n^2} - 1)}{\overset{1-q}{\mathbb{M}_n}}$, where $M_0$ is simple with a unique trace. $p_l$ is
full in $M$. \\ 

\par

$Proof\ of\ Theorem\ 2.6:$ \\ 

\par

Now to prove Theorem 2.6 we start with  

\begin{equation*}
(\mathfrak{A},\phi )=(\underset{\alpha_0}{\overset{p_0}{A_0}} \oplus \underset{\alpha_1}{\overset{p_1}{\mathbb{M}_{n_1}}} \oplus ...
\oplus \underset{\alpha_k}{\overset{p_k}{\mathbb{M}_{n_k}}})*(\underset{\beta_0}{\overset{q_0}{B_0}} \oplus 
\underset{\beta_1}{\overset{q_1}{\mathbb{M}_{m_1}}} \oplus ... \oplus \underset{\beta_l}{\overset{q_l}{\mathbb{M}_{m_l}}}),
\end{equation*}
where $A_0$ and $B_0$ have unital, diffuse abelian $C^*$-subalgebras, supported on their units 
(we allow $\alpha_0 = 0$ or/and $\beta_0 =
0$). The case where $n_1 = ... = n_k = m_1 = ... = m_l = 1$ is treated in Theorem 2.5. The case where 
$\alpha_0 = 0$, $k = 1$, and $n_k >
1$ was treated in Proposition 4.4. So we can suppose without loss of generality that $n_k \geq 2$ and either 
$k > 1$ or $\alpha_0 > 0$ or
both. To prove that the conclusions of Theorem 2.6 takes place in this case we will use induction on 
$\card \{ i | n_i \geq 2 \} + \card \{ j | m_j \geq 2 \}$, having Theorem 2.5 ($\card \{ i | n_i \geq 2 \} + 
\card \{ j | m_j \geq 2 \} = 0$)
as first step of the induction. We look at 

\begin{equation*}
(\mathfrak{B},\phi|_\mathfrak{B})=(\underset{\alpha_0}{\overset{p_0}{A_0}} \oplus \underset{\alpha_1}{\overset{p_1}{\mathbb{M}_{n_1}}} 
\oplus ... \oplus \underset{\alpha_{k-1}}{\overset{p_{k-1}}{\mathbb{M}_{n_{k-1}}}} \oplus \underset{\alpha_k}{\overset{p_k}{\mathbb{C}}}) 
* (\underset{\beta_0}{\overset{q_0}{B_0}} \oplus 
\underset{\beta_1}{\overset{q_1}{\mathbb{M}_{m_1}}} \oplus ... \oplus \underset{\beta_l}{\overset{q_l}{\mathbb{M}_{m_l}}}) \subset
(\mathfrak{A},\phi).
\end{equation*}

We suppose that Theorem 2.6 is true for $(\mathfrak{B},\phi|_\mathfrak{B})$ and we will prove it for $(\mathfrak{A},\phi )$. This will be
the induction step and will prove Theorem 2.6. 
\par
Denote $L_0^{\mathfrak{A}} \overset{def}{=} \{ (i,j)| \frac{\alpha_i}{n_i^2} + \frac{\beta_j}{m_j^2} = 1 \}$, $L_0^{\mathfrak{B}}
\overset{def}{=} \{ (i,j)| i \leq k-1$ and $\frac{\alpha_i}{n_i^2} + \frac{\beta_j}{m_j^2} = 1 \} \cup \{ (k,j) | \frac{\alpha_k}{1^2} + 
\frac{\beta_j}{m_j^2} = 1 \}$ and similarly $L_+^{\mathfrak{A}} \overset{def}{=} \{ (i,j)| \frac{\alpha_i}{n_i^2} + \frac{\beta_j}{m_j^2}
> 1 \}$, and $L_+^{\mathfrak{B}} \overset{def}{=} \{ (i,j)| i \leq k-1$ and $\frac{\alpha_i}{n_i^2} + \frac{\beta_j}{m_j^2} > 1 \} \cup \{ (k,j) | \frac{\alpha_k}{1^2} + 
\frac{\beta_j}{m_j^2} > 1 \}$. Clearly $L_0^{\mathfrak{A}} \cap \{ 1 \leq i \leq k-1 \} = L_0^{\mathfrak{B}} \cap \{ 1 \leq i \leq k-1 \}$
and similarly $L_+^{\mathfrak{A}} \cap \{ 1 \leq i \leq k-1 \} = L_+^{\mathfrak{B}} \cap \{ 1 \leq i \leq k-1 \}$. Let
$N_{\mathfrak{A}}(i,j) = max(n_i, m_j)$ and let $N_{\mathfrak{B}}(i,j) = N_{\mathfrak{A}}(i,j), 1 \leq i \leq k-1$, and 
$N_{\mathfrak{B}}(k,j) = m_j$.

By assumption 

\begin{equation*}
\mathfrak{B}= \underset{\delta}{\overset{g}{\mathfrak{B}_0}} \oplus \underset{(i,j)\in L_+^{\mathfrak{B}}}{\bigoplus}
\underset{\delta_{ij}}{\overset{g_{ij}}{\mathbb{M}_{N_{\mathfrak{B}}(i,j)}}}. 
\end{equation*}

We want to show that 

\begin{equation}
\mathfrak{A} = \underset{\gamma}{\overset{f}{\mathfrak{A}_0}} \oplus \underset{(i,j)\in L_+^{\mathfrak{A}}}{\bigoplus}
\underset{\gamma_{ij}}{\overset{f_{ij}}{\mathbb{M}_{N_{\mathfrak{A}}(i,j)}}}. 
\end{equation}

We can represent $\mathfrak{A}$ as the span of $p_k \mathfrak{A} p_k$, $p_k \mathfrak{A} p_k \mathfrak{B} (1-p_k)$, $(1-p_k) \mathfrak{B} p_k
\mathfrak{A} p_k$, $(1-p_k) \mathfrak{B} p_k \mathfrak{A} p_k \mathfrak{B} (1-p_k)$, and $(1-p_k) \mathfrak{B} (1-p_k)$. 
From the fact that $g_{kj} \leq p_k$ and $g_{ij} \leq 1-p_k, \forall 1 \leq i \leq k-1$ we see that $p_k \mathfrak{B} (1-p_k) = p_k
\mathfrak{B}_0 (1-p_k)$, $(1-p_k) \mathfrak{B} p_k = (1-p_k) \mathfrak{B}_0 p_k$, and $(1-p_k) \mathfrak{B} (1-p_k) = 
(1-p_k) \mathfrak{B}_0 (1-p_k) \oplus \underset{i \neq k}{\underset{(i,j) \in L_+^{\mathfrak{B}}}{\bigoplus}} \mathbb{M}_{N(i,j)}$. All this tells us that we can represent $\mathfrak{A}$ as the span of $p_k \mathfrak{A} p_k$, $p_k \mathfrak{A} p_k \mathfrak{B}_0 (1-p_k)$, $(1-p_k) \mathfrak{B}_0 p_k
\mathfrak{A} p_k$, $(1-p_k) \mathfrak{B}_0 p_k \mathfrak{A} p_k \mathfrak{B}_0 (1-p_k)$, $ (1-p_k) \mathfrak{B}_0 (1-p_k)$, and $\underset{i \neq k}{\underset{(i,j)\in L_+^{\mathfrak{B}}}{\bigoplus}} \underset{\delta_{ij}}{\overset{g_{ij}}{\mathbb{M}_{N(i,j)}}}$.
\par
In order to show that $\mathfrak{A}$ has the form (9), we need to look at $p_k \mathfrak{A} p_k$. From Proposition 4.1 we have 
$$p_k \mathfrak{A} p_k \cong (p_k \mathfrak{B} p_k) * \mathbb{M}_{n_k} \cong 
(\underset{\frac{\delta}{\alpha_k}}{\overset{g}{p_k \mathfrak{B}_0 p_k}} \oplus \underset{(k,j)\in L_+^{\mathfrak{B}}}{\bigoplus} 
\underset{\frac{\delta_{kj}}{\alpha_k}}{\overset{g_{kj}}{\mathbb{M}_{N(k,j)}}}) * \mathbb{M}_{n_k}.$$ 
Since by assumption $p_k \mathfrak{B}_0 p_k$ has a unital, diffuse abelian $C^*$-subalgebra, supported on $1_{p_k \mathfrak{B}_0 p_k}$, we can use Proposition 4.4 to determine the form of $p_k \mathfrak{A} p_k$. 
\par
Thus $p_k \mathfrak{A} p_k$: 
\par
(i) Is simple with a unique trace if whenever for all $1 \leq r \leq l$ with $N(k,r) = 1$ we have 
$\frac{\delta_{kr}}{\alpha_k} < 1 - \frac{1}{n_k^2}$. 
\par
(ii) Is an extension $0 \rightarrow I \rightarrow p_k \mathfrak{A} p_k \rightarrow \mathbb{M}_{n_k} 
\rightarrow 0$ if $\exists 1 \leq r \leq l$, with $N(k,r) = 1$, and 
$\frac{\delta_{kr}}{\alpha_k} = 1 - \frac{1}{n_k^2}$. 
Moreover $I$ is simple with a unique trace and has no central projections. 
\par
(iii) Has the form $p_k \mathfrak{A} p_k = I \oplus \underset{n_k^2(\frac{\delta_{kr}}{\alpha_k} - 1 + 
\frac{1}{n_k^2})}{\mathbb{M}_{n_k}}$, where $I$ is unital, simple with a unique trace whenever 
$\exists 1 \leq r \leq l$ with $N(k,r) = 1$, and $\frac{\delta_{kr}}{\alpha_k} > 1 - \frac{1}{n_k^2}$.
\par
By assumption $\delta_{ij} = N(i,j)^2 (\frac{\alpha_i}{n_i^2} +\frac{\beta_j}{m_j^2} - 1)$, 
so when $r$ satisfies the conditions of case (iii) above, then $m_r = 1$ and 
$n_k^2(\frac{\delta_{kr}}{\alpha_k} - 1 + \frac{1}{n_k^2}) = n_k^2(\frac{\alpha_k + \beta_r - 1}{\alpha_k} + 
\frac{1}{n_k^2} -1) = \frac{n_k^2}{\alpha_k}(\frac{\alpha_k}{n_k^2} + \frac{\beta_r}{1^2} - 1)$, 
just what we needed to show. Defining $\mathfrak{A}_0 \overset{def}{=} (1-(\underset{(i,j) 
\in L_+^{\mathfrak{A}}}{\oplus} f_{ij})) \mathfrak{A} (1-(\underset{(i,j) \in L_+^{\mathfrak{A}}}{\oplus} 
f_{ij}))$, we see that $\mathfrak{A}$ has the form (9).
\par
We need to study $\mathfrak{A}_0$ now. Since clearly $g \leq f$, we see that 
$\mathfrak{A} p_k \mathfrak{B}_0 = \mathfrak{A} p_k g \mathfrak{B}_0 = \mathfrak{A} g p_k \mathfrak{B}_0 = 
\mathfrak{A}_0 p_k \mathfrak{B}_0$ and similarly $\mathfrak{A} p_k \mathfrak{B}_0 = \mathfrak{A}_0 p_k 
\mathfrak{B}_0$. From this and from what we proved above follows that:

\begin{gather}
\mathfrak{A}_0 \text{ is the span of } p_k \mathfrak{A}_0 p_k,\ (1-p_k) \mathfrak{B}_0 p_k \mathfrak{A}_0 p_k, \\ \notag
p_k \mathfrak{A}_0 p_k \mathfrak{B}_0 (1-p_k),\ (1-p_k) \mathfrak{B}_0 p_k \mathfrak{A}_0 p_k \mathfrak{B}_0 (1-p_k), \text{ and } (1-p_k) \mathfrak{B}_0 (1-p_k). 
\end{gather}

We need to show that for each of the projections $f p_s$, $0 \leq s \leq k$ and $f q_t$, $1 \leq t \leq l$, we have a unital, diffuse 
abelian $C^*$-subalgebra of $\mathfrak{A}_0$, supported on it. The ones, supported on $f p_s$, $1 \leq s \leq k-1$ come from 
$(1-p_k) \mathfrak{B}_0 (1-p_k)$ by the induction hypothesis. The one with unit $f p_k$ comes from the representation 
$p_k \mathfrak{A} p_k \cong (p_k \mathfrak{B} p_k) * \mathbb{M}_{n_k}$ and Proposition 4.4. For $1 \leq s \leq l$ we have 

\begin{gather}
q_s \mathfrak{A} q_s \cong \underset{\frac{\gamma}{\beta_s}}{\overset{f q_s}{{q_s \mathfrak{A}_0 q_s}}} 
\oplus \underset{1 \leq i \leq k-1}{\underset{(i,s) \in L_+^{\mathfrak{A}}}{\bigoplus}} 
\underset{\frac{\gamma_{is}}{\beta_s}}{\overset{f_{is}}{\mathbb{M}_{N_{\mathfrak{A}}(i,s)}}} \oplus 
\underset{\frac{\gamma_{ks}}{\beta_s}}{\overset{f_{ks}}{\mathbb{M}_{N_{\mathfrak{A}}(k,s)}}}
\end{gather}

and

\begin{gather}
q_s \mathfrak{B} q_s \cong  \underset{\frac{\delta}{\beta_s}}{\overset{g q_s}{q_s \mathfrak{B}_0 q_s}} \oplus 
\underset{1 \leq i \leq k-1}{\underset{(i,s) \in L_+^{\mathfrak{B}}}{\bigoplus}} 
\underset{\frac{\delta_{is}}{\beta_s}}{\overset{g_{is}}{\mathbb{M}_{N_{\mathfrak{B}}(i,s)}}} \oplus 
\underset{\frac{\delta_{ks}}{\beta_s}}{\overset{g_{ks}}{\mathbb{M}_{N_{\mathfrak{B}}(k,s)}}}. 
\end{gather}

From what we showed above follows that for $1 \leq i \leq k-1$ we have $\gamma_{is} = \delta_{is}$ and 
$f_{is} = g_{is}$. If $(k,s) \notin L_+^{\mathfrak{B}}$, (or $\alpha_k < 1 - \frac{\beta_s}{m_s^2}$), 
then $(k,s) \notin L_+^{\mathfrak{A}}$ and by (11) and (12) we see that $gq_s = fq_s$ and so in $\mathfrak{A}_0$ we have a unital, diffuse abelian
$C^*$-subalgebra with unit $gq_s = fq_s$, which comes from $\mathfrak{B}_0$. If $(k,s) \in L_+^{\mathfrak{B}}$, then $gq_s \lvertneqq fq_s$
and since we have a unital, diffuse abelian $C^*$-subalgebra of $\mathfrak{A}_0$, supported on $gq_s$, comming from $\mathfrak{B}_0$, we
need only to find a unital, diffuse abelian $C^*$-subalgebra of $\mathfrak{A}_0$, supported on $fq_s - gq_s$ 
and its direct sum with the one supported on $gq_s$ will be a unital, diffuse abelian $C^*$-subalgebra of 
$\mathfrak{A}_0$, supported on $fq_s$. But from the form (11) and (12)
it is clear that $fq_s - gq_s \leq g_{ks}$, since from (11) and (12) $(f_{1s} + ... + f_{(k-1)s}) q_s \mathfrak{A} q_s (f_{1s} + ... + f_{(k-1)s}) = (g_{1s} + ... + g_{(k-1)s}) q_s \mathfrak{B} q_s (g_{1s} + ... + g_{(k-1)s})$. It is also clear then that $ fq_s - gq_s = f g_{ks} \leq p_k$, since $gq_s \perp g_{ks}$. We look for this $C^*$-subalgebra in 
$$p_k \mathfrak{A} p_k = \underset{\frac{\gamma}{\alpha_k}}{\overset{fp_k}{ p_k \mathfrak{A}_0 p_k }} \oplus 
\underset{(k,j)\in L_+^{\mathfrak{A}}}{\bigoplus} \underset{\frac{\gamma_{kj}}{\alpha_k}}{\overset{f_{kj}}
{\mathbb{M}_{N_{\mathfrak{A}}(k,j)}}} \cong (p_k \mathfrak{B} p_k) * \mathbb{M}_{n_k},$$ 
$$ \cong (\underset{\frac{\delta}{\alpha_k}}{\overset{g}{p_k \mathfrak{B}_0 p_k}} \oplus \underset{(k,j)\in L_+^{\mathfrak{B}}}{\bigoplus} 
\underset{\frac{\delta_{kj}}{\alpha_k}}{\overset{g_{kj}}{\mathbb{M}_{N_{\mathfrak{B}}(k,j)}}}) * 
\mathbb{M}_{n_k}.$$ 
Proposition 4.4 gives us a unital,
diffuse abelian $C^*$-subalgebra of $p_k \mathfrak{A}_0 p_k$, supported on $(f p_k) g_{ks} = f g_{ks} = 
fq_s -gq_s$. This proves that we have a unital, diffuse abelian $C^*$-subalgebra of $\mathfrak{A}_0$, 
supported on $fq_s$. 
\par
Now we have to study the ideal structure of $\mathfrak{A}_0$, knowing by the induction hypothesis, the form 
of $\mathfrak{B}$. We will
use the "span representation" of $\mathfrak{A}_0$ (10).
\par
For each $(i,j) \in L_0^{\mathfrak{B}}$ we know the existance of $*$-homomorphisms 
$\pi_{(i,j)}^{\mathfrak{B}_0} : \mathfrak{B}_0 \rightarrow \mathbb{M}_{N_{\mathfrak{B}}(i,j)}$. 
For $i \neq k$ we can write those as $\pi_{(i,j)}^{\mathfrak{B}_0} : \mathfrak{B}_0 \rightarrow 
\mathbb{M}_{N_{\mathfrak{A}}(i,j)}$ and since the support of $\pi_{(i,j)}^{\mathfrak{B}_0}$ is contained in 
$(1-p_k)$, using (10), we can extend linearly $\pi_{(i,j)}^{\mathfrak{B}_0}$ to $\pi_{(i,j)}^{\mathfrak{A}_0} : \mathfrak{A}_0 \rightarrow \mathbb{M}_{N_{\mathfrak{A}}(i,j)}$, by defining it to be zero on $p_k \mathfrak{A}_0 p_k$, $(1-p_k) \mathfrak{B}_0 p_k \mathfrak{A}_0 p_k$,  $p_k \mathfrak{A}_0 p_k \mathfrak{B}_0 (1-p_k)$, and $(1-p_k) \mathfrak{B}_0 p_k \mathfrak{A}_0 p_k \mathfrak{B}_0 (1-p_k)$. Clearly $\pi_{(i,j)}^{\mathfrak{A}_0}$ is a $*$-homomorphism also. 
\par
By the induction hypothesis we know that $g p_k$ is full in $\underset{i \neq k}{\underset{(i,j) 
\in L_0^{\mathfrak{B}}}{\bigcap}} \ker(\pi_{(i,j)}^{\mathfrak{B}_0}) \subset \mathfrak{B}_0$ and by (10), 
and the way we extended $\pi_{(i,j)}^{\mathfrak{B}_0}$, we see that $f p_k$ is full in 
$\underset{i \neq k}{\underset{(i,j) \in L_0^{\mathfrak{A}}}{\bigcap}} \ker(\pi_{(i,j)}^{\mathfrak{A}_0}) 
\subset \mathfrak{A}_0$. Then $p_k \mathfrak{A}_0 p_k$ is full and hereditary in 
$\underset{i \neq k}{\underset{(i,j) \in L_0^{\mathfrak{A}}}{\bigcap}} \ker(\pi_{(i,j)}^{\mathfrak{A}_0})$, 
so by the Rieffel correspondence from \cite{R82}, we have that $p_k \mathfrak{A}_0 p_k$ and 
$\underset{i \neq k}{\underset{(i,j) \in L_0^{\mathfrak{A}}}{\bigcap}} \ker(\pi_{(i,j)}^{\mathfrak{A}_0})$ 
have the same ideal structure.
\par
Above we saw that 

\begin{gather}
p_k \mathfrak{A} p_k = \underset{\frac{\gamma}{\alpha_k}}{\overset{fp_k}{ p_k \mathfrak{A}_0 p_k }} \oplus 
\underset{(k,j)\in L_+^{\mathfrak{A}}}{\bigoplus} \underset{\frac{\gamma_{kj}}{\alpha_k}}{\overset{f_{kj}}
{\mathbb{M}_{N_{\mathfrak{A}}(k,j)}}} \cong (p_k \mathfrak{B} p_k) * \mathbb{M}_{n_k} \cong \\ \notag
\cong (\underset{\frac{\delta}{\alpha_k}}{\overset{gp_k}{p_k \mathfrak{B}_0 p_k}} \oplus \underset{(k,j)\in L_+^{\mathfrak{B}}}{\bigoplus} 
\underset{\frac{\delta_{kj}}{\alpha_k}}{\overset{g_{kj}}{\mathbb{M}_{N_{\mathfrak{B}}(k,j)}}}) * \mathbb{M}_{n_k}.
\end{gather}

From Proposition 4.4 follows that $p_k \mathfrak{A}_0 p_k$ is not simple if and only if 
$\exists 1 \leq s \leq m$, such that $(k,s) \in L_+^{\mathfrak{B}}, m_s = 1$ with 
$\frac{\delta_{ks}}{\alpha_k} = 1-\frac{1}{n_k^2}$, where $\delta_{ks} = \alpha_k + \beta_s -1$. This means 
that $\frac{\alpha_k + \beta_s -1}{\alpha_k} = 1 - \frac{1}{n_k^2}$, which is equivalent to 
$\frac{\beta_s}{1^2} + \frac{\alpha_k}{n_k^2} = 1$, so this implies $(k,s) \in L_0^{\mathfrak{A}}$. If this 
is the case (13), together with Proposition 4.4 gives us a $*$-homomorphism 
$\pi'_{(k,s)} : p_k \mathfrak{A}_0 p_k \rightarrow \mathbb{M}_{n_k}$, such that 
$\ker(\pi'_{(k,s)}) \subset p_k \mathfrak{A}_0 p_k$ is simple with a unique trace. 
Using (10) we extend $\pi'_{(k,s)}$ linearly to a linear map $\pi_{(k,s)}^{\mathfrak{A}_0} : 
\mathfrak{A}_0 \rightarrow \mathbb{M}_{n_k}$, by defining $\pi_{(k,s)}^{\mathfrak{A}_0}$ to be zero on 
$(1-p_k) \mathfrak{B}_0 p_k \mathfrak{A}_0 p_k$, $p_k \mathfrak{A}_0 p_k \mathfrak{B}_0 (1-p_k)$, $(1-p_k) 
\mathfrak{B}_0 p_k \mathfrak{A}_0 p_k \mathfrak{B}_0 (1-p_k)$, and $(1-p_k) \mathfrak{B}_0 (1-p_k)$. 
Similarly as before, $\pi_{(k,s)}^{\mathfrak{A}_0}$ turns out to be a $*$-homomorphism. 
By the Rieffel correspondence of the ideals of $p_k \mathfrak{A}_0 p_k$ and 
$\underset{i \neq k}{\underset{(i,j) \in L_0^{\mathfrak{A}}}{\bigcap}} \ker(\pi_{(i,j)}^{\mathfrak{A}_0})$, 
it is easy to see that the simple ideal $\ker(\pi'_{(k,s)}) \subset p_k \mathfrak{A}_0 p_k$ corresponds to 
the ideal $\underset{(i,j) \in L_0^{\mathfrak{A}}}{\bigcap} \ker(\pi_{(i,j)}^{\mathfrak{A}_0}) \subset 
\underset{i \neq k}{\underset{(i,j) \in L_0^{\mathfrak{A}}}{\bigcap}} \ker(\pi_{(i,j)}^{\mathfrak{A}_0})$, so 
$\underset{(i,j) \in L_0^{\mathfrak{A}}}{\bigcap} \ker(\pi_{(i,j)}^{\mathfrak{A}_0})$ is simple. To see that 
$\underset{(i,j) \in L_0^{\mathfrak{A}}}{\bigcap} \ker(\pi_{(i,j)}^{\mathfrak{A}_0})$ has a unique trace we 
notice that from the construction of $\pi_{(i,j)}^{\mathfrak{A}_0}$ we have $\ker(\pi'_{(k,s)}) = p_k 
\ker(\pi_{(k,s)}^{\mathfrak{A}_0}) p_k = p_k \underset{(i,j) \in L_0^{\mathfrak{A}}}{\bigcap} 
\ker(\pi_{(i,j)}^{\mathfrak{A}_0}) p_k$ (the last equality is true because $p_k \mathfrak{A}_0 p_k \subset \
underset{i \neq k}{\underset{(i,j) \in L_0^{\mathfrak{A}}}{\bigcap}} \ker(\pi_{(i,j)}^{\mathfrak{A}_0})$). 
Now we argue similarly as in the proof of Proposition 4.3, using the fact that $\ker(\pi'_{(k,s)})$ has a 
unique trace: Suppose that $\rho$ is a trace on $\underset{(i,j) \in L_0^{\mathfrak{A}}}{\bigcap} 
\ker(\pi_{(i,j)}^{\mathfrak{A}_0})$. It is easy to see that $\Span \{ x p_k a p_k y | x, y, a \in 
\underset{(i,j) \in L_0^{\mathfrak{A}}}{\bigcap} \ker(\pi_{(i,j)}^{\mathfrak{A}_0}), a \geq 0 \}$ is dense in 
$\underset{(i,j) \in L_0^{\mathfrak{A}}}{\bigcap} \ker(\pi_{(i,j)}^{\mathfrak{A}_0})$, since 
$\ker(\pi'_{(k,s)})$ is full in $\underset{(i,j) \in L_0^{\mathfrak{A}}}{\bigcap} 
\ker(\pi_{(i,j)}^{\mathfrak{A}_0})$. Then since $p_k a p_k \geq 0$ we have $\rho(x p_k a p_k y) = 
\rho((p_k a p_k) y x) = \rho((p_k a p_k)^{1/2} y x (p_k a p_k)^{1/2})$ and since $(p_k a p_k)^{1/2} y x 
(p_k a p_k)^{1/2}$ is supported on $p_k$, it follows that $(p_k a p_k)^{1/2} y x (p_k a p_k)^{1/2} \in p_k 
\underset{(i,j) \in L_0^{\mathfrak{A}}}{\bigcap} \ker(\pi_{(i,j)}^{\mathfrak{A}_0}) p_k = \ker(\pi'_{(k,s)})$, 
so $\rho$ is uniquely determined by $\rho|_{\ker(\pi'_{(k,s)})}$ and hence $\underset{(i,j) 
\in L_0^{\mathfrak{A}}}{\bigcap} \ker(\pi_{(i,j)}^{\mathfrak{A}_0})$ has a unique trace.
\par
If $\nexists 1 \leq s \leq m$ with $(k,s) \in L_0^{\mathfrak{A}}$ it follows from what we said above, that 
$p_k \mathfrak{A}_0 p_k$ is simple with a unique trace. But since $p_k \mathfrak{A}_0 p_k$ is full and 
hereditary in $\underset{i \neq k}{\underset{(i,j) \in L_0^{\mathfrak{A}}}{\bigcap}} 
\ker(\pi_{(i,j)}^{\mathfrak{A}_0}) = \underset{(i,j) \in L_0^{\mathfrak{A}}}{\bigcap} 
\ker(\pi_{(i,j)}^{\mathfrak{A}_0})$ it follows  that $\underset{(i,j) \in L_0^{\mathfrak{A}}}{\bigcap} 
\ker(\pi_{(i,j)}^{\mathfrak{A}_0})$ is simple with a unique trace in this case too.
\par
We showed already that $f p_k$ is full in $\underset{i \neq k}{\underset{(i,j) 
\in L_0^{\mathfrak{A}}}{\bigcap}} \ker(\pi_{(i,j)}^{\mathfrak{A}_0})$. 
Now let $1 \leq r \leq k-1$. We need to show that $f p_r$ is full in 
$\underset{i \neq r}{\underset{(i,j) \in L_0^{\mathfrak{A}}}{\bigcap}} \ker(\pi_{(i,j)}^{\mathfrak{A}_0})$. 
From (11) and (12) follows that $f-g \leq p_k$. So $f p_r = g p_r$ for all $1 \leq r \leq k-1$. 
From the way we constructed $\pi_{(i,j)}^{\mathfrak{A}_0}$ is clear that $f p_r \in 
\underset{i \neq r}{\underset{(i,j) \in L_0^{\mathfrak{A}}}{\bigcap}} \ker(\pi_{(i,j)}^{\mathfrak{A}_0})$. 
It is also true that $f p_r \notin  \ker(\pi_{(r,j)}^{\mathfrak{A}_0})$ for any $1 \leq j \leq l$. So the 
smallest ideal of $\mathfrak{A}_0$, that contains $f p_r$, is $\underset{i \neq r}{\underset{(i,j) 
\in L_0^{\mathfrak{A}}}{\bigcap}} \ker(\pi_{(i,j)}^{\mathfrak{A}_0})$, meaning that we must have $\langle f 
p_r \rangle_{\mathfrak{A}_0} = \underset{i \neq r}{\underset{(i,j) \in L_0^{\mathfrak{A}}}{\bigcap}} 
\ker(\pi_{(i,j)}^{\mathfrak{A}_0})$.
\par
Finally, we need to show that for all $1 \leq s \leq l$ we have that $f q_s$ is full in 
$\underset{j \neq s}{\underset{(i,j) \in L_0^{\mathfrak{A}}}{\bigcap}} \ker(\pi_{(i,j)}^{\mathfrak{A}_0})$. 
Let $(i,j) \in L_0^{\mathfrak{A}}$ with $i \neq k$, $j \neq s$. 
Since $g q_s \in \ker(\pi_{(i,j)}^{\mathfrak{B}})$ and since $(f-g)q_s \leq p_k$, the way we extended 
$\pi_{(i,j)}^{\mathfrak{B}}$ to $\pi_{(i,j)}^{\mathfrak{A}}$ shows that 
$f q_s \in \ker(\pi_{(i,j)}^{\mathfrak{B}})$. Let $(i,s) \in L_0^{\mathfrak{A}}$ and $i \neq k$. 
Then we know that $g q_s \notin \ker(\pi_{(i,j)}^{\mathfrak{B}})$, which implies 
$f q_s \notin \ker(\pi_{(i,j)}^{\mathfrak{A}})$. Suppose $(k,s) \in L_0^{\mathfrak{A}}$. 
Then $m_s = 1$ and (13), Proposition 4.4,  and the way we extended $\pi'_{(k,s)}$ to 
$\pi_{(k,s)}^{\mathfrak{A}_0}$ show, that $f g_{ks} = fq_s - gq_s$ is full in $p_k \mathfrak{A}_0 p_k$, 
meaning that $fq_s -gq_s$, and consequently $fq_s$, is not contained in $\ker(\pi_{(k,s)}^{\mathfrak{A}_0})$. 
Finally let $j \neq s$, and suppose $(k,j) \in L_0^{\mathfrak{A}}$. This means that 
$(k,j) \in L_+^{\mathfrak{B}}$ and also that the trace of $q_j$ is so big, that 
$(i,s) \notin L_+^{\mathfrak{B}}$ and $(i,s) \notin L_0^{\mathfrak{B}}$ for any $1 \leq i \leq k$. 
Then (12) shows that $q_s \leq g$. The way we defined $\pi_{(k,j)}^{\mathfrak{A}_0}$ using (13) and 
Proposition 4.4 shows us that $\mathfrak{B}_0 \subset \ker(\pi_{(k,j)}^{\mathfrak{A}_0})$ in this case. 
This shows $q_s = g q_s = fq_s \in \ker(\pi_{(k,j)}^{\mathfrak{A}_0})$. All this tells us that the smallest 
ideal of $\mathfrak{A}_0$, containing $fq_s$, is 
$\underset{j \neq s}{\underset{(i,j) \in L_0^{\mathfrak{A}}}{\bigcap}} \ker(\pi_{(i,j)}^{\mathfrak{A}_0})$, 
and therefore $\langle fq_s \rangle_{\mathfrak{A}_0} = \underset{j \neq s}{\underset{(i,j) 
\in L_0^{\mathfrak{A}}}{\bigcap}} \ker(\pi_{(i,j)}^{\mathfrak{A}_0})$.
\par
This concludes the proof of Theorem 2.6.

\qed

{\em Acknowledgements.} I would like to thank Ken Dykema, my advisor, for the many helpful conversations I 
had with him, for the moral support and for reading the first version of this paper. I would also like to thank Ron Douglas and Roger Smith for some 
discussions.


\begin{thebibliography}{99}
\bibitem{APT73} C.\ Akeman,\ G.\ Pedersen,\ J.\ Tomiyama,
{\em Multipliers of $C^*$-Algebras, }
J.\ Func.\ Analysis, {\bf 13}, 1973, 277-301.
\bibitem{ABH91} J.\ Anderson\ B.\ Blackadar\ U.\ Haagerup,
{\em Minimal Projections in the Reduced Group $C^*$-algebra of $\mathbb{Z}_n*\mathbb{Z}_m$, }
J.\ Operator\ Theory, {\bf 26}, 1991, 3-23.
\bibitem{A82} D.\ Avitzour,
{\em Free Products of $C^*$-Algebras, }
Trans.\ Amer.\ Math.\ Soc. {\bf 271}, 1982, 423-435.
\bibitem{B93} E.\ B$\acute{e}$dos,
{\em On the Uniqueness of the Trace on Some Simple $C^*$-Algebras, }
J.\ Oper.\ Theory, {\bf 30}, 1993, 149-160.
\bibitem{C79} M.\ D.\ Choi,
{\em A Simple $C^*$-Algebra Generated by Two Finite-Order Unitaries, }
Canad.\ J.\ Math.\ {\bf 31},\ (1979),\ No.4,\ 867-880.
\bibitem{D93} K.\ Dykema,
{\em Free Products of Hyperfinite von Neumann Algebras and Free Dimension, }
Duke\ Math.\ J. {\bf 69}, No.1, 1993, 97-119. 
\bibitem{D94} K.\ Dykema,
{\em Interpolated Free Group Factors, }
Pacific\ J.\ Math. {\bf 163}, (1994), 123-134.
\bibitem{D98} K.\ Dykema,
{\em Faithfulness of free product states, }
J.\ Funct.\ Anal. {\bf 154},\ (1998),\ No.\ 2.
\bibitem{D99} K.\ Dykema,
{\em Simplicity and Stable Rank of Some Free Product $C^*$-Algebras,}
Trans.\ Amer.\ Math.\ Soc. {\bf 351}, No. 1, 1999, 1-40.
\bibitem{D99LN} K.\ Dykema,
{\em Free Probability Theory and Operator Algebras, }
Seoul National University GARC lecture notes, in preparation.
\bibitem{DR98} K.\ Dykema,\ M.\ R\o rdam,
{\em Projections in Free Product $C^*$-algebras, }
GAFA, Vol. {\bf 8}, (1998), 1-16.
\bibitem{DHR97} K.\ Dykema,\ U.\ Haagerup,\ M. R{\o}rdam,
{\em The Stable Rank of Some Free Product $C^*$-algebras, }
Duke\ Math.\ Journal, {\bf 90} No. 1, 1997, 95-121.
\bibitem{KR86} R.\ Kadison,\ J.\ Ringrose,
{\em Fundamentals of the Theory of Operator Algebras, }
Academic\ Press,\ Boston-New\ York,\ 1986.
\bibitem{Ka69} R.\ Kallman,
{\em A Generalization of Free Action, }
Duke\ Math.\ J.\ {\bf 36},\ 1969,\ 781-789.
\bibitem{K81} A.\ Kishimoto,
{\em Outer Automorphisms and Reduced Crossed Products of Simple $C^*$-Algebras,}
Commun.\ Math.\ Phys.\ {\bf 81}, 1981, 429-435.
\bibitem{L81} R.\ Longo,
{\em A Remark on Crossed Product of $C^*$-Algebras, }
J.\ London\ Math.\ Soc.\ {\bf 23}, (1981), 531-533.
\bibitem{O75} D.\ Olesen,
{\em Inner $*$-Automotrphisms of Simple $C^*$-Algebras,}
Commun.\ Math.\ Phys.\ {\bf 44}, 1975, 175-190.
\bibitem{OP78} D.\ Olesen,\ G.\ Pedersen,
{\em Application of the Connes Spectrum to $C^*$-Dynamical Systems,}
J.\ Func.\ An.\ {\bf 30}, 1978, 179-197.
\bibitem{PS79} W.\ Paschke,\ N.\ Salinas,
{\em $C^*$-Algebras Associated with Free Products of Groups, }
Pacific\ J.\ Math.\ {\bf 82},\ (1979),\ No.1,\ 211-221.
\bibitem{P75} R.\ Powers,
{\em Simplicity of the $C^*$-Algebra, Associated with the Free Group
on Two Generators, }
Duke\ Math.\ J.\ {\bf 42},\ (1975),\ 151-156.
\bibitem{Ra94} F.\ R\u{a}dulescu,
{\em Random Matrices, Amalgamated Free Products and Subfactors of the von Neumann Algebra of a Free Group, of Noninteger Index, }
Invent.\ Math. {\bf 115}, (1994), 347-389.
\bibitem{R82} M.\ Rieffel,
{\em Morita Equivalence for Operator Algebras, }
Proc.\ Symp.\ Pure\ Math. {\bf 38},\ (1982),\ 285-298.
\bibitem{S81} S.\ Str\u{a}til\u{a},
{\em Modular Theory in Operator Algebras, }
Editura Academiei-Abacus Press,\ 1981.
\bibitem{T79} M.\ Takesaki, 
{\em Theory of Operator Algebras I, }
Springer-Verlag,\ New\ York-Heidelberg-Berlin,\ 1979.
\bibitem{V85} D.\ Voiculescu,
{\em Symmetries of Some Reduced Free Product $C^*$-Algebras, }
Operator Algebras and Their Connection with Topology and Ergodic Theory, Lecture Notes in Mathematics, Volume {\bf 1132}, Springer-Verlag, 1985, 556-588. 
\end{thebibliography}
\end{document}